\input amstex
\documentstyle{amsppt}
\TagsOnRight
\pagewidth{5.5 in}
\pageheight{8 in}
\magnification=1200
\NoBlackBoxes
 
\def\conj#1{\overline{#1}}
\def\inner#1#2#3{ \langle #1,\, #2 \rangle_{_{#3}}}
\def\Dinner#1#2{ \langle #1,\, #2 \rangle_{D}}
\def\Dpinner#1#2{ \langle #1,\, #2 \rangle_{D^\perp}}

\def\vartau{\tau\hskip2pt\!\!\!\!\iota\hskip1pt}
\def\ft#1{\widehat{#1}}
	
\def\BB{1}
\def\FB{2}
\def\MD{3}
\def\JD{4}
\def\EE{5}
\def\EL{6}
\def\EG{7}
\def\KS{8}
\def\HLa{9}
\def\HLb{10}
\def\HLc{11}
\def\CP{12}
\def\MRc{13}
\def\TU{14}
\def\SWa{15}
\def\SWb{16}
\def\SWc{17}
\def\SWd{18}
\def\WW{19}

\topmatter

\title
The AF structure of non commutative \\ toroidal $\Bbb Z/4\Bbb Z$ orbifolds
\endtitle

\rightheadtext{AF structure of non commutative toroidal orbifolds}

\author
S.~Walters
\endauthor

\affil
\sevenrm The University of Northern British Columbia
\\ \vskip5pt
{\sevenit Dedicated to my Mother on her seventieth}
\endaffil

\date
\sevenrm July 25, 2002
\enddate

\address
Department of Mathematics and Computer Science, The University
of Northern British Columbia, Prince George, B.C.  V2N 4Z9  CANADA
\endaddress

\email
walters\@hilbert.unbc.ca \ or \ walters\@unbc.ca  \hfill
\break\indent{\it Home page}: http://hilbert.unbc.ca/walters
\endemail

\thanks
Research partly supported by NSERC grant OGP0169928
\hfill {\sevenrm (\TeX File: orbifoldLANL.tex)}
\endthanks



\abstract
For any irrational $\theta$ and rational number $p/q$ such that $q|q\theta-p|<1$, 
a projection $e$ of trace $q|q\theta-p|$ is constructed in the the irrational rotation 
algebra $A_\theta$ that is invariant under the Fourier transform.  (The latter is
the order four automorphism $U\mapsto V \mapsto U^{-1}$, where $U,V$ are the canonical 
unitaries generating $A_\theta$.)
Further, the projection $e$ is approximately central, the cut down algebra 
$eA_\theta e$ contains a Fourier invariant $q\times q$
matrix algebra whose unit is $e$, and the cut downs $eUe, eVe$ are 
approximately inside the matrix algebra. (In particular, there are
Fourier invariant projections of trace $k|q\theta-p|$ for $k=1,\dots,q$.)
It is also shown that for all $\theta$ the crossed product 
$A_\theta \rtimes \Bbb Z_4$ satisfies the Universal Coefficient Theorem. 
($\Bbb Z_4 := \Bbb Z/4\Bbb Z$.)
As a consequence, using the Classification Theorem of G.~Elliott and G.~Gong for 
AH-algebras, a theorem of M.~Rieffel, and by recent results of H.~Lin, we show that 
$A_\theta \rtimes \Bbb Z_4$ is an AF-algebra for all irrational $\theta$ in a dense $G_\delta$.
\endabstract

\endtopmatter


\document

\newpage

\subhead \S1. INTRODUCTION \endsubhead

For each irrational number $\theta$, the irrational rotation C*-algebra $A_\theta$ is
generated by unitaries $U,V$ enjoying the Heisenberg commutation relation  $VU=\lambda
UV$, where $\lambda=e^{2\pi i\theta}$.   The Fourier transform on $A_\theta$ is the
order four automorphism $\sigma$ (or its inverse) of $A_\theta$ given by
$\sigma(U)=V$ and $\sigma(V)=U^{-1}$. 
It is a non commutative analogue of the Fourier transform of classical analysis.
The resulting toroidal $\Bbb Z/4\Bbb Z$ orbifold $A_\theta^\sigma$ is the mathematical
physicist's name for the fixed point subalgebra of $A_\theta$ under $\sigma$.
Such algebras have found their way in recent work of Konechny and Schwarz 
[\KS].  One of the problems of George Elliott is to determine whether this fixed point 
subalgebra, or the C*-crossed product $A_\theta \rtimes \Bbb Z_4$, is an AF-algebra.
In this paper we solve this problem in the affirmative for a dense $G_\delta$ set of 
irrationals $\theta$.  In order to do this we prove the following result.

\proclaim{Theorem 1.1}
Let $\theta$ be any irrational number in $(0,1)$.
Let $\tfrac pq$ be a rational number in reduced form (with $p\ge0,\, q\ge1$) such that
$|\theta - \tfrac pq| < \frac1{q^2}$.  (For example, this is satisfied by any 
convergent of $\theta$.)
Then there is a projection $e$ in $A_\theta$ of trace $q|q\theta-p|$
enjoying the following properties:
\itemitem{(1)} $e$ is invariant under the Fourier transform,
\itemitem{(2)} $e$ is approximately central, 
\itemitem{(3)} $e$ is the unit of a Fourier invariant $q\times q$ matrix algebra $\Cal M_q$
contained in $eA_\theta e$, 
\itemitem{(4)} the cut downs $eUe, eVe$ are approximately in $\Cal M_q$.
\endproclaim

Previous constructions of Fourier invariant projections placed some specialized 
assumptions on $\tfrac pq$ and on the bound $\tfrac1{q^2}$.  
In [\FB], Proposition 2.1, Boca proved
that if $p$ satisfies a quadratic residue condition modulo $q$ and $|\theta-\tfrac pq|
<\tfrac{0.948}{q^2}$, then there is a Fourier invariant projection of trace $|q\theta-p|$.
Also, in [\SWc], Corollary 1.2, the author assumed $q$ is a sum of two squares in 
constructing a projection orthogonal to its Fourier transform, and in
[\SWd] a result similar to Theorem 1.1 is obtained, for a dense $G_\delta$ 
set of irrationals $\theta$, in which one has a Fourier invariant direct sum of two matrix
algebras interchanged by the Fourier transform, but in which it is assumed that 
$q|q\theta-p|$ goes to zero as $q\to\infty$.
The projection constructed in this paper is a generalized Rieffel projection, i.e.~it 
has the C*-valued inner product form $\Dinner\xi\xi$ for a suitable Schwartz function 
$\xi$, a technique that appears to have been introduced by Elliott and Q.~Lin [\EL].

Once Theorem 1.1 is established, we obtain the following (see Section 9).

\proclaim{Theorem 1.2}
For a dense $G_\delta$ set of irrationals $\theta$ the C*-algebra
$A_\theta \rtimes \Bbb Z_4$ is an AF-algebra.
\endproclaim

This result stands on the shoulders of three major theorems: the Elliott-Gong 
Classification Theorem [\EG] for real rank zero, stable rank one, AH-algebras;
on Rieffel's Theorem 2.15 [\MRc] on strong Morita equivalence of twisted 
group C*-algebras associated to lattices in locally compact Abelian groups;
and on H.~Lin's recent theorem that every 
unital separable simple nuclear tracially AF C*-algebra that satisfies the Universal 
Coefficient Theorem (UCT) is isomorphic to an AH-algebra with slow dimension growth 
[\HLa] (as well as on Lin's paper [\HLb]).
A reading of the first few paragraphs of Section 9 below conveys more clearly how these 
results are applied in the present paper to obtain Theorem 1.2.  
Also shown in Section 9 is that the C*-algebra
$A_\theta \rtimes \Bbb Z_4$ satisfies UCT for all $\theta$ (Theorem 9.4).
We note that in fact the very same proof shows that the crossed products 
$A_\theta \rtimes_\alpha \Bbb Z_6$ and $A_\theta \rtimes_{\alpha^2}\Bbb Z_3$, where 
$\alpha$ is the canonical order 6 automorphism, $\alpha(U)=V,\ \alpha(V)=U^{-1}V$,
satisfy UCT for any $\theta$. 
The proof uses theorems of Tu [\TU] and Dadarlat [\MD].

As a consequence of Theorem 1.1, we have:

\proclaim{Corollary 1.3} Let $\theta$ be any irrational number, $0<\theta<1$, and let
$\tfrac pq, \tfrac{p'}{q'}$ be any two consecutive convergents of $\theta$ such that
$\tfrac pq < \theta < \tfrac{p'}{q'}$.  If $q'<q$, then there is a Fourier invariant 
projection in $A_\theta$ of trace $q'(q\theta-p)$ which is also the unit of a Fourier 
invariant $q'\times q'$ matrix algebra.  If $q'>q$, then there is a Fourier invariant 
projection in $A_\theta$ of trace $q(p'-q'\theta)$ which is also the unit of a Fourier 
invariant $q\times q$ matrix algebra.
\endproclaim

The projection of Corollary 1.3 is a subprojection of the projection of Theorem 1.1,
and can probably be chosen to be approximately central (but we have not checked this).
In fact, it is not hard to see that such a subprojection also satisfies (3) and (4) 
of Theorem 1.1 (for appropriate $q$).  (See Section 4.)

The plan of the paper is as follows:

\itemitem{$\bullet$} 
{\it Section 2}: we choose the appropriate Rieffel framework [\MRc] on which our 
results are based, as well as declare the Schwartz function $f$ that the construction 
depends on.

\itemitem{$\bullet$}
{\it Section 3}: the ``discrete part'' of the function $f$ is designed so that the
diophantine systems that later emerge will have a unique solution. 

\itemitem{$\bullet$}
{\it Section 4}: we establish the Fourier invariance property of the projection and for
the associated Fourier invariant matrix algebra, of which it is a unit.  

\itemitem{$\bullet$}
{\it Section 5}: an inner product calculation is carried out for Sections 6 and 8.

\itemitem{$\bullet$}
{\it Section 6}: the invertibility of the inner product $\Dpinner ff$ is established. This
yields the existence of the projection, and which, together with the invariance 
result, will give parts (1) and (3) of Theorem 1.1.  

\itemitem{$\bullet$}
{\it Section 7}: the projection is shown to be approximately central, giving (2) 
of Theorem 1.1.

\itemitem{$\bullet$}
{\it Section 8}: the projection is shown to be a point-----thus, the cut downs of the 
canonical unitaries are approximately in the invariant $q\times q$ matrix algebra-----giving
(4) of Theorem 1.1. 

\itemitem{$\bullet$}
{\it Section 9}: we show that the crossed product $A_\theta\rtimes\Bbb Z_4$ satisfies
UCT for each $\theta$ (Theorem 9.4), that it is tracially AF for a specifiable dense
$G_\delta$ set of irrationals $\theta$ (Theorem 9.1), and eventually that it is an 
AF-algebra for a dense $G_\delta$ (Theorem 9.3). 

The main Theta functions to be used in this paper are
$$
\vartheta_2(z,t) = \sum_n e^{\pi it(n+\frac 1 2)^2} e^{i2z(n+\tfrac12)}, \ \ \ 
\vartheta_3(z,t) = \sum_n e^{\pi itn^2} e^{i2zn}, \ \ \ 
\vartheta_4(z,t) = \sum_n (-1)^n e^{\pi itn^2} e^{i2zn}, 
$$
for $z,t\in\Bbb C$ and $\roman{Im}(t)>0$, where all summations range over the integers
$\Bbb Z$. (For a classic treatment see [\WW], whose definitions for $\vartheta_j$ we 
follow.)

We shall freely adopt the convention $e(t) := e^{2\pi i t}$.  For integers $n,q$ we write
$\delta_q^n = 1$ if $q|n$ and $\delta_q^n = 0$ otherwise.
Thus, $\sum_{j=0}^{q-1} e(\tfrac {nj}q) = q\delta_q^n$ for all integers $n$.

The author is indebted to George Elliott for mentioning the Fourier transform problem 
to him in 1994 and for many helpful e-mail exchanges since then.
The author is also indebted to N.~Christopher Phillips for valuable discussions related to 
classification and for a helpful suggestion in proving Theorem 9.4.  This research is partly
supported by a grant from NSERC.

\remark{Postscript}
After completing this paper we received a preliminary report from Chris Phillips [\CP] 
that he has been able to prove, among other things, that the C*-algebra 
$A_\theta \rtimes \Bbb Z_4$ is tracially AF for all irrational $\theta$ (generalizing our
Theorem 9.1). 
\endremark

\subhead \S2. RIEFFEL'S FRAMEWORK \endsubhead

The construction is based on Rieffel's Theorem 2.15 in [\MRc] stating that the Schwartz
space $\Cal S(M)$ on a locally compact Abelian group $M$ is an equivalence bimodule
with the C*-algebras $C^*(D,\frak h)$ and $C^*(D^\perp,\conj{\frak h})$, acting on the
left and the right, respectively.  We shall write these algebras simply as $C^*(D)$ and
$C^*(D^\perp)$. Here, $D$ is a lattice in $G=M\times \widehat M$, and
$\frak h$ the Heisenberg cocycle on $G$.  Further, $\pi: G \to \Cal L(L^2(M))$ is the 
Heisenberg representation given by $[\pi_{(m,s)}f](n) = \langle n,s\rangle f(n+m)$.
The C*-algebra generated by the unitaries $\pi_x$, for $x\in D$, is $C^*(D)$,
while $C^*(D^\perp)$ is the opposite algebra of the C*-algebra generated
by the unitaries $\pi_y^*$, for $y\in D^\perp$.  ($D^\perp = \{y\in G: 
\frak h(x,y) \conj{\frak h(y,x)} = 1,\ \ \forall x\in D\}$.)  We shall assume that 
$M$ is canonically isomorphic to its dual $\widehat M$ so that there is an order four map
$R:G\to G$ given by $R(u;v)=(-v;u)$.  A good summary of this situation can be found in 
Boca's paper [\FB] (or in [\SWc]).  In [\FB] it is shown that on each of these 
C*-algebras there is a canonical order four automorphism, corresponding
to the Fourier transform, given by $\sigma(\pi_x) = \conj{\frak h(x,x)} \pi_{Rx}$
and $\sigma'(\pi_y) = \conj{\frak h(y,y)} \pi_{Ry}$, for $x\in D,\ y\in D^\perp$.  One has
$$
\sigma(\Dinner fg) = \Dinner{\ft f}{\ft g}, \qquad
\sigma'(\Dpinner fg) = \Dpinner{\ft f}{\ft g}
$$
where $\ft f$ is the Fourier transform of $f$ (in the sense of harmonic analysis).
If $\vartau$ and $\vartau'$ are the canonical normalized traces on 
$C^*(D)$ and $C^*(D^\perp)$, respectively, then one has	
$$
\vartau( \Dinner f g ) = |G/D|\, \vartau' (\Dpinner g f)
$$
where $|G/D|$ is the Haar-Plancheral measure of a fundamental domain for $D$ in $G$.

Let $\theta$ be an irrational number with $0<\theta<1$.  Let $p/q$ be any
rational number, where $p\ge0, q\ge1$ are relatively prime, satisfying 
$$
\left| \theta - \frac pq \right| \ < \ \frac1{q^2}.	\tag2.1
$$
In particular, any convergent $p/q$ of $\theta$ satisfies this inequality.
We will construct invariant projections of trace $k|q\theta-p|$ for $k=1,\dots,q$.
Observe that with no loss of generality we can assume that $\tfrac pq < \theta$.
For if $\theta < \tfrac pq$, then $\tfrac{p'}q < 1-\theta$, where $p'=q-p$,
so that $q(q(1-\theta)-p')=q(p-q\theta)$.  In addition, the canonical isomorphism
$A_\theta \to A_{1-\theta}$ is compatible with the Fourier transform on each algebra,
so that their Fourier invariant projections are in one-to-one correspondence.

\bigpagebreak

Henceforth we will assume that $0< \tfrac pq < \theta$ satisfies (2.1).
By Lagrange's Theorem, one can write $p$ as a sum of four squares: 
$p = p_1^2 + p_2^2 + p_3^2 +p_4^2$, where $p_j\ge0$ are integers.  
Let $M=\Bbb R \times \Bbb Z_q \times \Bbb Z_q$ and consider the lattice $D$ in 
$G=M\times \widehat M$ with basis
$$
D:\ \ \bmatrix \varepsilon_1 \\ \varepsilon_2 \endbmatrix 
=
\bmatrix
\alpha & [p_1] & [p_2] & 0 & [p_3] & [p_4] \\
0 & [-p_3] & [-p_4] & \alpha & [p_1] & [p_2]
\endbmatrix,
$$
where $\alpha = (\theta -\tfrac pq)^{1/2}$, and $[n]:=[n]_q$ is the mod $q$ class of 
$n$.  It is clear that $D$ is invariant under the ``Fourier'' map $R$ on $G$, after 
$M$ and $\widehat M$ have been identified.  Since a fundamental domain for $D$ is
$[0,\alpha) \times \Bbb Z_q \times \Bbb Z_q \times 
[0,\alpha) \times \Bbb Z_q \times \Bbb Z_q$,
the covolume of $D$ is $|G/D| = \alpha^2q^2 = q^2\theta-pq < 1$ by (2.1).
The C*-algebra $C^*(D)$ is generated by the canonical unitaries
$\pi_{\varepsilon_1}, \pi_{\varepsilon_2}$ whose commutation relation is
$$
\pi_{\varepsilon_1} \pi_{\varepsilon_2} \pi_{\varepsilon_1}^* \pi_{\varepsilon_2}^*
= \frak h (\varepsilon_1, \varepsilon_2) \conj{\frak h(\varepsilon_2, \varepsilon_1)}
= e(\alpha^2+\tfrac{p_1^2 + p_2^2 + p_3^2 + p_4^2}q) = e(\theta) = \lambda.
$$
The unitaries 
$U_1 = \mu_0^{1/2} \pi_{\varepsilon_1},\ U_2 = \mu_0^{-1/2} \pi_{\varepsilon_2}$,
where $\mu_0 := e(\tfrac{p_1p_3+p_2p_4}q)$, enjoy the condition
$U_1U_2=\lambda U_2U_1$ and the Fourier transform can be checked to satisfy
$\sigma(U_1)=U_2, \sigma(U_2)=U_1^*$.

Throughout the paper, we shall let
$$
\beta := \frac1{q\alpha} \ > \ 1.
$$
The complementary lattice $D^\perp$ can easily be checked to be generated by the 
(dependent) elements
$$
D^\perp: \ \ 
\bmatrix \delta_1' \\ \\ \delta_2' \\ \\ \delta_3' \\ \\ \delta_4' \\ \\ \delta_5' 
\\ \\ \delta_6' \endbmatrix
=
\bmatrix
q\beta & 0 & 0 & 0 & 0 & 0 \\ \\
0 & 0 & 0 & q\beta & 0 & 0 \\ \\
-\beta p_1 & [1] & 0 & \beta p_3 & 0 & 0  \\ \\
-\beta p_2 & 0 & [1] & \beta p_4 & 0 & 0  \\ \\
-\beta p_3 & 0 & 0 & -\beta p_1 & [1] & 0  \\ \\
-\beta p_4 & 0 & 0 & -\beta p_2 & 0 & [1]  \\
\endbmatrix.
$$
Let $c,d$ be integers such that $cp+dq=1$.
A little reflection will show that a basis for $D^\perp$ is given by
$$
D^\perp: \ \
\bmatrix \delta_1 \\ \\ \delta_2 \\ \\ \delta_3 \\ \\ \delta_4 \endbmatrix
=
\bmatrix
\beta & [-cp_1] & [-cp_2] & 0 & [-cp_3] & [-cp_4] \\ \\
0 & [cp_3] & [cp_4] & \beta & [-cp_1] & [-cp_2] \\ \\
0 & [p_2] & [-p_1] & 0 & [-p_4] & [p_3] \\ \\
0 & [p_4] & [-p_3] & 0 & [p_2] & [-p_1] 
\endbmatrix.
$$
Indeed, this follows from the following relations and their inverses
$$
\xalignat 2
\delta_1 &= d\delta_1' -cp_1 \delta_3' - cp_2 \delta_4' - cp_3 \delta_5' - cp_4 \delta_6',
&\qquad 
\delta_1' &= q\delta_1
\\
\delta_2 &= d\delta_2' +cp_3 \delta_3' + cp_4 \delta_4' - cp_1 \delta_5' - cp_2 \delta_6',
&\qquad 
\delta_2' &= q\delta_2
\\
\delta_3 &= p_2\delta_3' - p_1\delta_4' - p_4\delta_5' + p_3\delta_6',	
&\qquad 
\delta_3' &= -p_1 \delta_1 + p_3 \delta_2 + cp_2\delta_3 + cp_4\delta_4
\\
\delta_4 &= p_4\delta_3' - p_3\delta_4' + p_2\delta_5' - p_1\delta_6',
&\qquad 
\delta_4' &= -p_2 \delta_1 + p_4 \delta_2 - cp_1\delta_3 - cp_3\delta_4
\\
 &
&\qquad 
\delta_5' &= -p_3 \delta_1 - p_1 \delta_2 - cp_4\delta_3 + cp_2\delta_4
\\
 &	
&\qquad 
\delta_6' &= -p_4 \delta_1 - p_2 \delta_2 + cp_3\delta_3 - cp_1\delta_4
\\
\endxalignat
$$
which are easily checked using $p=\sum_jp_j^2$ and $cp+dq=1$.  To check that
$\{\delta_j\}$ is an independent set, it clearly suffices to check that 
$\delta_3,\delta_4$ are independent as order $q$ elements.  However, this is immediate.
Writing
$$
\sum_{j=1}^4 n_j\delta_j = (\beta n_1, [m_1], [m_2];\ \beta n_2, [m_3], [m_4])
$$
where
$$
\aligned
m_1 &= -cp_1n_1 + cp_3n_2 + p_2n_3 + p_4n_4 \\
m_2 &= -cp_2n_1 + cp_4n_2 - p_1n_3 - p_3n_4 \\
m_3 &= -cp_3n_1 - cp_1n_2 - p_4n_3 + p_2n_4 \\
m_4 &= -cp_4n_1 - cp_2n_2 + p_3n_3 - p_1n_4 \\
\endaligned
\tag2.2
$$
the $D^\perp$ inner product coefficients become (for any Schwartz functions $f,g$ on $M$)
$$
\align
&\Dpinner f g (\Sigma_j n_j\delta_j)
\\
&\ \ \ =
\Dpinner f g (\beta n_1, [m_1], [m_2];\ \beta n_2, [m_3], [m_4])
\\
&\ \ \ =
\int_{\Bbb R \times \Bbb Z_q \times \Bbb Z_q}
\!\!\!\!\conj{ f(x,[n],[m]) }  g(x+\beta n_1, [n+m_1], [m+m_2]) 
e(\beta xn_2 + \tfrac{nm_3+mm_4}q) \ dx d[n] d[m]
\\
&\ \ \ =
\frac1q \sum_{\Sb [n],[m]\in\Bbb Z_q\endSb}
e(\tfrac{nm_3+mm_4}q)
\int_{\Bbb R} \conj{ f(x,[n],[m]) } g(x+\beta n_1, [n+m_1], [m+m_2]) e(\beta n_2x) \,dx.
\endalign
$$
Hence
$$
\Dpinner fg = 
\sum_{n_1n_2n_3n_4} \Dpinner fg (\Sigma n_j\delta_j) \cdot \pi_{\Sigma n_j\delta_j}^*.
$$
Throughout the paper we shall consider the Schwartz function on $M$ defined by
$$
f(x,n,m) = \frac{2^{1/4}}{\sqrt q} \, h(x) \varphi(n,m), \qquad
\text{where} \qquad
\cases \varphi(n,m) = e(\tfrac1q[an^2+bnm+ \gamma m^2]) & \\
	h(x) = e^{-\pi x^2} &
\endcases
$$
and $a,b,\gamma$ are suitable integers to be chosen carefully in the next section.
One has
$$
\Dpinner ff (\Sigma n_j\delta_j)
=
\frac{\sqrt2}{q^2} \, \Omega \, H(n_1\beta,n_2\beta)
$$
where
$$
\Omega = 
\sum_{m,n=0}^{q-1} e(\tfrac{nm_3+mm_4}q) \conj{\varphi(n,m)} \varphi(n+m_1,m+m_2)
$$
and
$$
H(s,t) = \int_{\Bbb R} \conj{h(x)} h(x+s) e(tx) \,dx \ = \
\frac1{\sqrt2} e(-\tfrac12st) e^{-\tfrac\pi2(s^2+t^2)}.
$$
In the following it will be useful to use the identity
$$
\pi_{n\delta} = \frak h(\delta,\delta)^{-n(n-1)/2}\, \pi_{\delta}^n
$$
which holds for all vectors $\delta$ and integers $n$. Let 
$$
V_j := \lambda_{jj}^{-(q+1)/2} \pi_{\delta_j}^* = \lambda_{jj}^{-(q-1)/2} \pi_{-\delta_j}
$$  
for $j=1,2,3,4$, where $\lambda_{jk} = \frak h(\delta_j,\delta_k)$.  From
$\pi_{\delta_j} \pi_{\delta_k} = \lambda_{jk} \bar\lambda_{kj} \, \pi_{\delta_k} \pi_{\delta_j}$
one has
$$
V_j V_k = \lambda_{jk} \bar\lambda_{kj} \, V_k V_j
$$
and the commutation relations
$$
V_1V_2 = e(\theta') V_2V_1, \qquad V_3V_4 = e(\tfrac pq)\,V_4V_3, \qquad
V_jV_k=V_kV_j, \qquad V_k^q = I,
$$
for $j=1,2$ and $k=3,4$, where 
$\theta' = \beta^2 + \tfrac cq = \tfrac{c\theta+d}{q\theta-p}$ is in the
standard $\roman{GL}(2,\Bbb Z)$ orbit of $\theta$.  Hence $C^*(D^\perp)$ is isomorphic
to $M_q(A_{\theta'})$ (by universality and simplicity of the latter).  One can also
check that the Fourier transform $\sigma'$ on this C*-algebra is given by
$$
\sigma'(V_1) = V_2, \quad \sigma'(V_2) = V_1^*, \quad
\sigma'(V_3) = V_4, \quad \sigma'(V_4) = V_3^*.
$$
From the relations
$$
\frak h(u+v,w) = \frak h(u,w) \frak h(v,w), \qquad
\frak h(u,v+w) = \frak h(u,v) \frak h(u,w),
$$
and since $\lambda_{14}=\lambda_{23}=1$, one gets 
$$
\frak h(n_1\delta_1+n_2\delta_2,\, n_3\delta_3+n_4\delta_4) =
\lambda_{13}^{n_1n_3} \lambda_{24}^{n_2n_4}
$$
so that
$$
\align
\pi_{\Sigma n_j\delta_j} 
&= \conj{\frak h(n_1\delta_1+n_2\delta_2,\, n_3\delta_3+n_4\delta_4)}
\ \pi_{n_1\delta_1+n_2\delta_2} \pi_{n_3\delta_3+n_4\delta_4} 
\\
&= \lambda_{13}^{-n_1n_3} \lambda_{24}^{-n_2n_4} 
\ \conj{\frak h(n_1\delta_1, n_2\delta_2)} \pi_{n_1\delta_1} \pi_{n_2\delta_2}
\ \conj{\frak h(n_3\delta_3, n_4\delta_4)} \pi_{n_3\delta_3} \pi_{n_4\delta_4}
\\
&=
\lambda_{13}^{-n_1n_3} \lambda_{24}^{-n_2n_4} \lambda_{12}^{-n_1n_2} \lambda_{34}^{-n_3n_4}
\pi_{n_1\delta_1} \pi_{n_2\delta_2}\pi_{n_3\delta_3} \pi_{n_4\delta_4}.
\endalign
$$
Hence
$$
\pi_{\Sigma n_j\delta_j}^* \ = \ 
\Lambda_{n_1n_2n_3n_4}\,V_4^{n_4} V_3^{n_3} V_2^{n_2} V_1^{n_1}
$$
where
$$
\Lambda_{n_1n_2n_3n_4} := 
\lambda_{13}^{n_1n_3} \lambda_{24}^{n_2n_4}
\lambda_{12}^{n_1n_2} \lambda_{34}^{n_3n_4}
\lambda_{11}^{n_1(n_1+q)/2}\, \lambda_{22}^{n_2(n_2+q)/2}\,
\lambda_{33}^{n_3(n_3+q)/2}\, \lambda_{44}^{n_4(n_4+q)/2}.
\tag2.3
$$
Thus we have
$$
\Dpinner fg = 
\sum_{n_1n_2n_3n_4} \Dpinner fg (\Sigma n_j\delta_j) 
\, \Lambda_{n_1n_2n_3n_4}\,V_4^{n_4} V_3^{n_3} V_2^{n_2} V_1^{n_1}.
\tag2.4
$$

\bigpagebreak

\newpage

\subhead \S3. A DIOPHANTINE SYSTEM \endsubhead

With $p/q$ being any positive rational and $p=p_1^2+p_2^2+p_3^2+p_4^2$, and with $m_j$ 
as defined in (2.2), consider the diophantine system of congruences in $n_j$ ($j=1,2,3,4$)
$$
m_3+2am_1+bm_2 + u_3 \equiv 0, \qquad m_4+2\gamma m_2+bm_1 + u_4 \equiv 0
\tag3.1
$$
which will arise in Section 5, where (throughout the rest of the paper) we shall 
simply write ``$n\equiv m$'' for $n\equiv m \mod q$, and where $u_3,u_4$ are integers.
Inserting the values of $m_j$, the diophantine system (3.1) becomes
$$
\aligned
r_4 n_3 + r_2 n_4 &\equiv -u_3 + cs_3 n_1 + cs_1 n_2, \\
r_3 n_3 + r_1 n_4 &\equiv -u_4 + cs_4 n_1 + cs_2 n_2,
\endaligned
\tag3.2
$$
where 
$$
\xalignat 2
r_1 &= -p_1-2\gamma p_3 + bp_4,	&\qquad s_1 &= p_1-2ap_3-bp_4
\\
r_2 &= p_2+2ap_4-bp_3, 		&\qquad	s_2 &= p_2 - 2\gamma p_4 - bp_3   \tag3.3
\\
r_3 &= p_3- 2\gamma p_1 + bp_2, &\qquad	s_3 &= p_3 + 2a p_1 + bp_2 
\\
r_4 &= -p_4+ 2a p_2 - bp_1, 	&\qquad	s_4 &= p_4 + 2\gamma p_2 + bp_1.
\\
\endxalignat
$$
In this section we shall arrange for the determinant $\Delta = r_1r_4-r_2r_3$ of 
the system (3.2) to be relatively prime to $q$.  Once this is done, we can then 
solve (3.2) for $n_3,n_4$ and obtain
$$
n_3 = c_3 + a_1n_1 + a_2n_2, \qquad n_4 = c_4 + b_1n_1 + b_2n_2,
\tag3.4
$$
where
$$
\aligned
a_1 &= c \Delta'(r_1s_3-r_2s_4), \qquad a_2 = c\Delta'(r_1s_1-r_2s_2), 
\\
b_1 &= c \Delta'(r_4s_4-r_3s_3), \qquad b_2 = c\Delta'(r_4s_2-r_3s_1),
\\
c_3 &= \Delta'(-r_1u_3+r_2u_4), \qquad c_4 = \Delta'(-r_4u_4+r_3u_3),
\endaligned
\tag3.5
$$
and $\Delta'$ is an integer such that $\Delta'\Delta \equiv 1 \mod q$.
It is straightforward to check that one also has $\Delta = s_1s_4-s_2s_3$ by (3.3).

More explicity, the determinant is
$$
\split
\Delta &:= 
(-p_1-2\gamma p_3+bp_4)(-p_4+2ap_2-bp_1) - (p_2+2ap_4-bp_3)(p_3-2\gamma p_1+bp_2)
\\
&\ = bA + 2(\gamma-a)B + (1-b^2+4a\gamma)C
\endsplit
$$
where
$$
A = p_1^2 - p_2^2 + p_3^2 - p_4^2, \qquad B = p_1p_2+p_3p_4, \qquad
C = p_1p_4-p_2p_3.
$$
The following beautiful $ABC$ relation is easy to check:
$$
A^2 + 4B^2 + 4C^2 = p^2.
$$
From this it follows that $(A,B,C)$ divides $p$, and hence is relatively prime to $q$.

\bigpagebreak

\proclaim{Proposition 3.1} Let $p/q$ be any positive rational in reduced form and
write $p=p_1^2+p_2^2+p_3^2+p_4^2$ (by Lagrange's Theorem).  Then there are integers 
$a,b,\gamma$ such that $(\Delta,q)=1$.
\endproclaim

To prove this, we need the following two lemmas.	

\proclaim{Lemma 3.2}
Let $p,r,q$ be integers such that $(p,r,q)=1$.  Then there is an integer
$k$ such that $(p+kr,q)=1$.
\endproclaim

\demo{Proof}
Let $t_1,\dots,t_n$ be the distinct prime divisors of $q$ that divide some
element of the set $\{p+mr: m\in\Bbb Z\}$.  Then no $t_j$ divides $(p,r)$.
Let $k_j$ be an integer such that $t_j$ divides $p+k_jr$.  
If $n=1$, letting $k=k_1+1$, it is easy to see that $t_1$ does not divide $p+kr$.
(Since otherwise $t_1|r$ and hence $t_1|p$, a contradiction.)
Now suppose $n\ge2$.
Let $d_j = (k_1-k_j, t_j)$, for $j=2,\dots,n$ and put
$
k = k_1 + \frac{t_2}{d_2} \cdots \frac{t_n}{d_n}.
$
First we claim that $t_j \nmid (k-k_j)$ for $j=1,\dots,n$.  Suppose for some
$j$ one has $t_j | (k-k_j)$.  For $j=1$ this is clearly impossible.  Hence 
assume $j\ge2$.  Then $t_j$ divides
$k_1-k_j + \tfrac{t_2}{d_2} \cdots \tfrac{t_j}{d_j} \cdots \tfrac{t_n}{d_n}$
which is not possible whether $d_j=1$ or $t_j$, as is easily seen, hence the
claim.  Now if $t_j |(p+kr)$for some $j$, then $t_j | (k-k_j)r$ and hence
$t_j|r$.  But then $t_j|p$ and so $t_j|(p,r)$, a contradiction.  \qed
\enddemo

\proclaim{Lemma 3.3}
Let $X,Y,Z$ be integers such that $(X,Y,Z)=1$.
Then there is an integer $k$ such that $kX+(k^2-1)Y$ is relatively prime to $Z$ .
\endproclaim

\demo{Proof} Let $t_1,\dots,t_m$ be the prime divisors of $Z$ that do not divide $Y$ 
and $s_1,\dots,s_n$ the prime divisors of $(Z,Y)$.
Let $k=t_1\cdots t_m$ and consider $\beta := kX + k^2Y-Y$.  If some $t_j|\beta$, then 
$t_j|Y$, a contradiction.  If some $s_i|\beta$, then $s_i|kX$.  But as $s_i\nmid k$, 
$s_i$ must divide $X$ and hence also $(X,Y,Z)=1$, another contradiction.\qed
\enddemo

\demo{Proof of Proposition 3.1} The proof is divided according to the parity of $q$.

Suppose $q$ is odd.
Let $b=2ka$ and $\gamma=k^2a$ where $k,a$ are integers to be determined.  Then
one has $\Delta = 2ka A + 2a(k^2-1)B + C = 2a\beta_k + C$, where 
$\beta_k := kA+(k^2-1)B$.
Applying Lemma 3.3 with $X=A,\ Y=B,\ Z=(q,C)$, one obtains an integer $k$ such that 
$\beta_k$ is relatively prime to $(q,C)$.  In particular, the hypothesis of Lamma 3.2,
$(2\beta_k,C,q)=1$, is now satisfied (as $q$ is odd), hence there is an integer 
$a$ such that $\Delta = 2a\beta_k + C$ is relatively prime to $q$.

Suppose $q$ is even.  
In this case $A$ is odd since $p$ is odd using the $ABC$ relation above.
Let $b=1+2ka$ and $\gamma=k^2a$ where $k,a$ are integers to be determined.
Then one checks that
$\Delta = 2a\beta_k + A$ where $\beta_k := k(A-2C)+(k^2-1)B$.  Applying Lemma 3.3 with 
$X=A-2C,\ Y=B,\ Z=(q,A)$ (which are relatively prime since $A$ is odd), one obtains an 
integer $k$ such that $\beta_k$, and hence $2\beta_k$, is relatively prime to $(q,A)$.
Again Lemma 3.2 now yields an integer $a$ such that $\Delta = 2a\beta_k + A$ is 
relatively prime to $q$.  \qed
\enddemo

\bigpagebreak

Henceforth, and throughout the rest of the paper, $a,b,\gamma$ will be fixed integers
such that $\Delta$ is relatively prime to $q$.

\newpage

\subhead \S4. INVARIANCE \endsubhead

We begin by looking at the above groups/lattices obtained by removing the copy of 
$\Bbb R$.  Thus, $M_0 = \Bbb Z_q \times \Bbb Z_q$,\ $G_0 = M_0 \times \ft M_0 =
M_0 \times M_0$, \ $D_0$ the lattice in $G_0$, and $D_0^\perp$ its complement, with bases:
$$
D_0:\ \ \bmatrix \varepsilon_1 \\ \varepsilon_2 \endbmatrix
=
\bmatrix
p_1 & p_2 & p_3 & p_4 \\
-p_3 & -p_4 & p_1 & p_2
\endbmatrix,
\qquad
D_0^\perp: \ \
\bmatrix \delta_3 \\ \delta_4 \endbmatrix
=
\bmatrix
p_2 & -p_1 &  -p_4 & p_3 \\ 
p_4 & -p_3 &  p_2 & -p_1
\endbmatrix.
$$
The fact that $\{\varepsilon_1, \varepsilon_2,  \delta_3, \delta_4\}$ is a basis
for $G_0$ follows since the determinant
$$
\det \bmatrix
p_1 & p_2 & p_3 & p_4 \\
-p_3 & -p_4 & p_1 & p_2 \\
p_2 & -p_1 &  -p_4 & p_3 \\ 
p_4 & -p_3 &  p_2 & -p_1
\endbmatrix
= -(p_1^2+p_2^2+p_3^2+p_4^2)^2 = -p^2
$$
is relatively prime to $q$.  It is easily checked that $\delta_3,\delta_4\in D_0^\perp$,
and that if $n_1\varepsilon_1+n_2\varepsilon_2 \in D_0^\perp$, then $n_1 \equiv n_2 
\equiv 0 \mod q$, hence $D_0^\perp$ does indeed have $\{\delta_3,\delta_4\}$ as a basis.
Further, a fundamental domain for $D_0$ is $D_0^\perp$ itself because the restriction of
the canonical map $G_0 \to G_0/D_0$ to $D_0^\perp$ is clearly a bijection.  Since in terms 
of the aforementioned basis $D_0^\perp$ is bijective with 
$\{0\} \times \{0\} \times \Bbb Z_q \times \Bbb Z_q$, one obtains the covolume
$|G_0/D_0| = \tfrac1{\sqrt q} \tfrac1{\sqrt q} \sqrt q \sqrt q = 1$.
Now the C*-algebra $C^*(D_0) \cong M_q$ is generated by unitaries $u_1,u_2$ with
$u_1u_2 = e(\tfrac pq) u_2u_1$, and $C^*(D_0^\perp) \cong M_q$ is 
generated by unitaries $V_3,V_4$ of order $q$ enjoying $V_3V_4=e(\tfrac pq)V_4V_3$.
On each of these algebras we have respective
Fourier transforms $\sigma_0(\pi_x) = \conj{\frak h(x,x)} \pi_{R_0x}$ for $x\in D_0$,
and $\sigma_0'(\pi_y) = \conj{\frak h(y,y)} \pi_{R_0y}$ for $y\in D_0^\perp$, where
$R_0(u;v) = (-v;u)$ for $u,v\in M_0$.  They satisfy the properties 
$$
\sigma_0( \inner{\phi_1}{\phi_2}{D_0} ) = \inner{\ft\phi_1}{\ft\phi_2}{D_0}, \qquad
\sigma_0'( \inner{\phi_1}{\phi_2}{D_0^\perp} ) = \inner{\ft\phi_1}{\ft\phi_2}{D_0^\perp}.  
$$
If $\vartau_0, \vartau_0'$ are the canonical normalized traces on $C^*(D_0)$ and
$C^*(D_0^\perp)$, respectively, then 
$\vartau_0( \inner{\phi_1}{\phi_2}{D_0} ) = 
\vartau_0'( \inner {\phi_2}{\phi_1}{D_0^\perp} )$, for $\phi_j \in \Cal(M_0)$, since
we have $|G_0/D_0|=1$.  

\proclaim{Lemma 4.1} If $\varphi(n,m) = e(\tfrac1q[an^2+bnm+ \gamma m^2])$, 
with $a,b,\gamma$ chosen according to Proposition 3.1, then there is a unitary
$W_0$ in $C^*(V_3,V_4) \cong M_q(\Bbb C)$ such that 
$$
\ft\varphi = \varphi W_0.
$$
In fact, $W_0 = \inner\varphi{\ft\varphi}{D_0^\perp}$, after suitably normalizating 
$\varphi$ so that $\inner\varphi\varphi{D_0^\perp} = 1$.  Further, $\sigma_0'(W_0)=W_0^*$.
\endproclaim

\demo{Proof}
We will show that, after suitably normalizating $\varphi$,
one has $\inner\varphi\varphi{D_0^\perp} = 1$.  By the above trace relation it will
also follow that $\inner\varphi\varphi{D_0} = 1$ (being a projection of trace 1), and
hence applying the Fourier transforms one has 
$\inner{\ft\varphi}{\ft\varphi}{D_0^\perp} = 1$ and
$\inner{\ft\varphi}{\ft\varphi}{D_0} = 1$.  Therefore, one obtains
$\ft\varphi = \inner\varphi\varphi{D_0} \ft\varphi 
= \varphi \inner\varphi{\ft\varphi}{D_0^\perp}$.  Next one notes that $W_0 = 
\inner\varphi{\ft\varphi}{D_0^\perp}$ is a unitary in the matrix algebra $C^*(V_3,V_4)$
since 
$$
W_0 W_0^*  
= \inner\varphi{\ft\varphi}{D_0^\perp} \inner\varphi{\ft\varphi}{D_0^\perp}^*
= \inner\varphi{\ft\varphi}{D_0^\perp} \inner{\ft\varphi}\varphi{D_0^\perp}
= \inner\varphi{\ft\varphi \inner{\ft\varphi}\varphi{D_0^\perp}}{D_0^\perp} 
= \inner\varphi{ \inner{\ft\varphi}{\ft\varphi}{D_0} \varphi } {D_0^\perp} = 1.
$$
Therefore, we need only show  $\inner\varphi\varphi{D_0^\perp}$ is a scalar.

\noindent For conveneince, let
$$
M=p_2k+p_4\ell, \quad N=-p_1k-p_3\ell, \quad M'=-p_4k+p_2\ell, \quad N'=p_3k-p_1\ell.
$$
One has
$$
\align
&\inner \varphi\varphi{D_0^\perp}(k\delta_3+\ell\delta_4) 
=
\inner\varphi\varphi{D_0^\perp}(M,N;M',N')
\\
&=
\int_{M_0} \conj{\varphi(n,m)} \varphi(n+M,m+N)  e(\tfrac{nM'+mN'}q) d[n]d[m]
\\
&= \sum_{m,n=0}^{q-1}  e(\tfrac{nM'+mN'}q)  \conj{\varphi(n,m)} \varphi(n+M,m+N)
\\
&=
\sum_{m,n=0}^{q-1}  e(\tfrac{nM'+mN'}q) \conj{e(\tfrac1q[an^2+bnm+\gamma m^2])}
\\
&\hskip1in \cdot
 e(\tfrac1q[a(n+M)^2+b(n+M)(m+N)+\gamma(m+N)^2])
\\
&=
\sum_{m,n=0}^{q-1}  e(\tfrac{nM'+mN'}q) \conj{e(\tfrac1q[an^2+bnm+\gamma m^2])}
\\
&\ \ \ \ \ \ \ \ 
\cdot e(\tfrac1q[a(n^2+2Mn+M^2)+b(nm+Mm+Nn+MN)+\gamma(m^2+2Nm+N^2)])
\\
&=
\sum_{m,n=0}^{q-1}  e(\tfrac{nM'+mN'}q)
e(\tfrac1q[a(2Mn+M^2)+b(Mm+Nn+MN)+\gamma(2Nm+N^2)])
\\
&=
e(\tfrac1q[aM^2+bMN+\gamma N^2])
\sum_{m,n=0}^{q-1}  e(\tfrac{nM'+mN'}q)  e(\tfrac1q[2aMn+bMm+bNn+2\gamma Nm])
\\
&=
e(\tfrac1q[aM^2+bMN+\gamma N^2])
\sum_{m,n=0}^{q-1}  e(\tfrac nq [M'+2aM+bN])  e(\tfrac mq [N'+2\gamma N+bM])
\\
&=
q^2 e(\tfrac1q[aM^2+bMN+\gamma N^2]) \delta_q^{M'+2aM+bN} \delta_q^{N'+2\gamma N+bM}.
\endalign
$$
These coefficients are nonzero iff the diophantine system
$$
M'+2aM+bN \equiv 0 \mod q, \qquad N'+2\gamma N+bM \equiv 0 \mod q,
$$
holds, or when written out they become, in terms of (3.3),
$$
r_4k + r_2\ell \ \equiv \ 0 \mod q, \qquad
r_3k + r_1\ell \ \equiv \ 0 \mod q.
\tag4.1
$$
Its determinant is $\Delta=r_1r_4-r_2r_3$ of Section 3.  By Proposition 3.1, 
$\Delta$ is relatively prime to $q$, so it follows that (4.1) holds iff 
$k \equiv \ell \equiv 0 \mod q$.  This means
that $\inner \varphi\varphi{D_0^\perp}$ is a positive scalar, and hence the result.  \qed
\enddemo

\newpage

Since the unitaries $V_3,V_4$ of the matrix algebra act on the right of 
$f=h\varphi$ only by acting on $\varphi$ (since they do not affect $h$), which is
why we used the same notation for them as for the $V_3,V_4$ of Section 2,
the same will hold of the unitary $W_0$, and by Lemma 4.1 one obtains (as $\ft h = h$) 
$$
\ft f = h \ft\varphi = h (\varphi W_0) = (h\varphi)W_0 = fW_0.  \tag4.2
$$
In Section 6 it will be shown that $b^{-2} := \Dpinner ff$ is invertible.
Once this is done one sets $\xi := fb$ so that $\Dpinner\xi\xi = 1$ and 
$e := \Dinner\xi\xi$ is a projection in $C^*(D) \cong A_\theta$ of trace 
$q(q\theta-p)$.  One then has the isomorphism $\mu_\xi: eA_\theta e \to C^*(D^\perp)$
given by $\mu_\xi(x) = \Dpinner \xi{x\xi}$, with inverse 
$\mu_\xi^{-1}(y) = \Dinner {\xi y}\xi$.  One easily checks that 
$\mu_{\ft\xi} \sigma = \sigma' \mu_\xi$. From (4.2) one has 
$$
\sigma'( \Dpinner ff ) = \Dpinner{\ft f}{\ft f} = \Dpinner{fW_0}{fW_0} = 
W_0^* \Dpinner ff W_0.
$$
Thus, $\sigma'(b) = W_0^* bW_0$, and one gets 
$$
\ft\xi = \ft f \sigma'(b) = (fW_0)(W_0^* bW_0) = (fb)W_0 = \xi W_0.
$$
Therefore, we can deduce that $e$ is a Fourier invariant projection:
$$
\sigma(e) = \sigma(\Dinner\xi\xi) = \Dinner{\ft\xi}{\ft\xi} = \Dinner{\xi W_0}{\xi W_0}
= \Dinner\xi\xi = e,
\tag4.3
$$
where we used the adjoint property $\Dinner{\xi a}{\eta} = \Dinner{\xi}{\eta a^*}$
for equivalence bimodules.  Next we show that the inverse image of $C^*(V_3,V_4) 
\cong M_q$ under $\mu_\xi$ is an invariant $q\times q$ matrix algebra under the 
Fourier transform $\sigma$.  Note that
$$
\mu_{\ft\xi}(x) = \Dpinner{\ft\xi}{x\ft\xi} = \Dpinner{\xi W_0}{x\xi W_0}
= W_0^* \Dpinner \xi{x\xi} W_0 = W_0^* \mu_\xi(x) W_0 = \rho(\mu_\xi(x))
$$
where $\rho(\cdot) = W_0^*(\cdot)W_0$ is an automorphism of $C^*(D^\perp) \cong 
M_q(A_{\theta'})$ that is the identity on $A_{\theta'}$ and leaves the matrix
factor invariant.
Since we already have $\mu_{\ft\xi} \sigma = \sigma' \mu_\xi$, one
gets $\mu_\xi \sigma = (\rho^{-1}\sigma') \mu_\xi$, from which it
is immediate that the inverse image of $C^*(V_3,V_4) \cong M_q$ under $\mu_\xi$ is
invariant under $\sigma$ (as $C^*(V_3,V_4)$ is invariant under $\sigma'$ and $\rho$).
To summarize, we have obtained the following result.

\proclaim{Theorem 4.2} The projection $e=\Dinner\xi\xi$, where $\xi=fb$ and
$b=\Dpinner ff^{-1/2}$, is invariant under the Fourier transform $\sigma$, has
trace $q|q\theta-p|$, and is the unit of a Fourier invariant $q\times q$ matrix 
algebra contained in $A_\theta$.  In particular, there is a Fourier invariant 
subprojection $e'$ of $e$ of trace $k|q\theta-p|$ for $k=1,2,\dots,q-1$ that is the unit
of a Fourier invariant $k\times k$ matrix subalgebra of $e'A_\theta e'$.
\endproclaim

\proclaim{Corollary 4.3} Let $\theta$ be any irrational number, $0<\theta<1$, and let
$\tfrac pq, \tfrac{p'}{q'}$ be any two consecutive convergents of $\theta$ such that
$\tfrac pq < \theta < \tfrac{p'}{q'}$.  If $q'<q$, then there is a Fourier invariant 
projection in $A_\theta$ of trace $q'(q\theta-p)$ which is also the unit of a Fourier 
invariant $q'\times q'$ matrix algebra.  If $q'>q$, then there is a Fourier invariant 
projection in $A_\theta$ of trace $q(p'-q'\theta)$ which is also the unit of a Fourier 
invariant $q\times q$ matrix algebra.
\endproclaim

\subhead Primitive Form of the Projection \endsubhead
Now let us calculate the inner product $\Dinner ff$ where $f$ is
the Schwartz function of Section 2.  We will need this for Section 7.  First, We have
$$
\split
\Dinner fg &(m\varepsilon_1 + n \varepsilon_2) =
\Dinner fg (m\alpha, mp_1-np_3, mp_2-np_4;\ n\alpha, mp_3+np_1, mp_4+np_2)
\\
&=
\frac1q
\sum_{r,s=0}^{q-1} \int_{\Bbb R} f(x,r,s) \conj{ g(x+m\alpha, r+mp_1-np_3, s+mp_2-np_4) }
\\
&\hskip1in \cdot e(-xn\alpha) e(-\tfrac{r(mp_3+np_1)}q) \, e(-\tfrac{s(mp_4+np_2)}q) \ dx 
\endsplit
$$
and
$$
\Dinner fg = |G/D| \sum_{m,n} \Dinner fg (m\varepsilon_1 + n \varepsilon_2) 
\pi_{m\varepsilon_1 + n \varepsilon_2}
$$
but
$$
\pi_{m\varepsilon_1 + n \varepsilon_2} =
\conj{\frak h(m\varepsilon_1, n \varepsilon_2)} \pi_{m\varepsilon_1} \pi_{n \varepsilon_2}
=
e(-\alpha^2-\tfrac{p_1^2+p_2^2}q)^{mn} \pi_{m\varepsilon_1} \pi_{n \varepsilon_2}
$$
Since for any $\varepsilon$ one has
$$
\pi_{m\varepsilon} = \conj{\frak h(\varepsilon,\varepsilon)}^{m(m-1)/2}\pi_{\varepsilon}^m
$$
one gets
$$
\split
\pi_{m\varepsilon_1 + n \varepsilon_2} 
&=
e(-\alpha^2-\tfrac{p_1^2+p_2^2}q)^{mn} 
\conj{\frak h(\varepsilon_1,\varepsilon_1)}^{m(m-1)/2}\pi_{\varepsilon_1}^m
\conj{\frak h(\varepsilon_2,\varepsilon_2)}^{n(n-1)/2}\pi_{\varepsilon_2}^n
\\
&=
e(-\theta+\tfrac pq-\tfrac{p_1^2+p_2^2}q)^{mn}
\conj{\frak h(\varepsilon_1,\varepsilon_1)}^{m(m-1)/2} 
\conj{\frak h(\varepsilon_2,\varepsilon_2)}^{n(n-1)/2} \, \pi_{\varepsilon_1}^m 
\pi_{\varepsilon_2}^n
\endsplit
$$
Now
$$
\frak h(\varepsilon_1,\varepsilon_1) = e(\tfrac{p_1p_3+p_2p_4}q) \ =: \mu_0, \qquad
\frak h(\varepsilon_2,\varepsilon_2) = e(\tfrac{-p_1p_3-p_2p_4}q) = \conj \mu_0, 
$$
and since $U_1 = \mu_0^{1/2}\pi_{\varepsilon_1}$ and 
$U_2 = \mu_0^{-1/2}\pi_{\varepsilon_2}$ one obtains
$$
\pi_{m\varepsilon_1 + n \varepsilon_2} = 
e(-\theta+\tfrac{p_3^2+p_4^2}q)^{mn} \mu_0^{(n^2-m^2)/2} U_1^m U_2^n
=
\mu_{m,n} U_2^n U_1^m
$$
where
$$
\mu_{m,n} \ := \ e(\tfrac{p_3^2+p_4^2}q)^{mn}\ \mu_0^{(n^2-m^2)/2}.
$$
Since $|G/D|=1/\beta^2$, we have
$$
\Dinner fg = \frac1{\beta^2} \sum_{m,n} \Dinner fg (m\varepsilon_1 + n \varepsilon_2)
\mu_{m,n} \, U_2^n U_1^m.
$$
For convenience, let
$$
M=mp_1-np_3, \quad N=mp_2-np_4, \quad M'=mp_3+np_1, \quad N'=mp_4+np_2.
$$
Then
$$
\align
\Dinner ff &(m\varepsilon_1 + n \varepsilon_2)
\\
&= 
\frac{\sqrt2}{q^2} \sum_{r,s=0}^{q-1} e(-\tfrac{(rM'+sN')}q) \, \int h(x) \varphi(r,s) 
\conj{h(x+m\alpha) \varphi(r+M, s+N)} \, e(-xn\alpha) dx
\\
&=
\frac{\sqrt2}{q^2} \conj{H(m\alpha,n\alpha)} \sum_{r,s=0}^{q-1} e(-\tfrac{(rM'+sN')}q) 
e(\tfrac1q[ar^2+brs+\gamma s^2]) 
\\
&\hskip1.5in \cdot e(-\tfrac1q[a(r+M)^2+b(r+M)(s+N)+\gamma(s+N)^2])
\\
&=
\frac{\sqrt2}{q^2}  \conj{H(m\alpha,n\alpha)} e(-\tfrac1q[aM^2+bMN+\gamma N^2]) 
\\
&\hskip1in \cdot
\sum_{r=0}^{q-1} e(-\tfrac rq[M'+2aM+bN])
\sum_{s=0}^{q-1} e(-\tfrac sq[N'+2\gamma N+bM]) 
\\
&=
\sqrt2\, \conj{H(m\alpha,n\alpha)} e(-\tfrac1q[aM^2+bMN+\gamma N^2])
\cdot \delta_q^{M'+2aM+bN} \delta_q^{N'+2\gamma N+bM}.
\endalign
$$
This inner product coefficient is nonzero exactly when the diophantine system of congruences
(in terms of the notation in (3.3))
$$
\aligned
M'+2aM+bN &= s_1n + s_3m \equiv 0 \mod q
\\
N'+2\gamma N+bM &= s_2n + s_4m\equiv 0 \mod q
\endaligned
$$
has a solution for $m,n$. Its determinant is $\Delta = s_1s_4-s_2s_3$
which is the same as that of the diophantine system considered in
Section 3, and which is relatively prime to $q$.
This shows that $\delta_q^{M'+2aM+bN} \delta_q^{N'+2\gamma N+bM}=1$ iff 
$m \equiv n \equiv 0 \mod q$ (and in which case $M \equiv N \equiv 0 \mod q$). Thus,
$\Dinner ff (qm\varepsilon_1 + qn \varepsilon_2) = \sqrt2 \, \conj{H(qm\alpha,qn\alpha)}$
and
$$
\Dinner ff = \frac{\sqrt2}{\beta^2}
\sum_{m,n} \conj{H(qm\alpha,qn\alpha)} \mu_{qm,qn} U_2^{qn} U_1^{qm}.
$$
where $\mu_{qm,qn} = (-1)^{q(p_1p_3+p_2p_4)(n^2-m^2)}$. Since we took 
$h(x)=e^{-\pi x^2}$ one has
$$
H(s,t) = \frac1{\sqrt2} e(-\tfrac12 st) e^{-\tfrac\pi2(s^2+t^2)}
$$
so that
$$
X := \Dinner ff =
\frac1{\beta^2}
\sum_{m,n} e(\tfrac12q^2\alpha^2mn) e^{-\tfrac\pi2q^2\alpha^2(m^2+n^2)} 
\mu_{qm,qn} U_2^{qn} U_1^{qm}.
$$

\bigpagebreak

Now we refer to this element as a {\it primitive form} of the projection $e$ because it is
positive and is supported on $e$.  That is, (i) $Xe=X$, and (ii) $ye=0$ if, and only if, 
$yX=0$ for any positive element $y$.  
To see (i), note that $e=\Dinner{fb}{fb} = \Dinner{fb^2}f$, and 
$$
Xe = \Dinner ff \Dinner{fb^2}f = \Dinner{\Dinner ff fb^2}f =
\Dinner{ f \Dpinner f{fb^2}} f = \Dinner ff = X
$$
by definition of $b$ (recall $b^{-2}:=\Dpinner ff$).
For (ii), if $ye=0$, then $yX=0$ by (i).  So assume $yX=0$.
Then $0=yXy = \Dinner{yf}{yf}$, so that $yf=0$.
Hence $yfb=0$ and so $0=\Dinner{yfb}{fb} = y \Dinner {fb}{fb} = ye$, as
needed.

\newpage

\subhead \S5. INNER PRODUCT CALCULATION \endsubhead

In this section we calculate the inner products $\Dpinner ff$ and $\Dpinner f{U_1f}$
in one swoop for Sections 6 and 8.  Since
$$
(\pi_{\varepsilon_1}f)(x,n,m) = e(\tfrac{p_3n+p_4m}q)\, f(x+\alpha,n+p_1,m+p_2)
$$
we look more generally at the function
$$
g(x,n,m) = e(\tfrac{t_3n+t_4m}q)\, f(x+t,n+t_1,m+t_2)
$$
where $t$ is real and $t_j$ are integers, which gives $f$ for $t=t_j=0$ and gives
$\pi_{\varepsilon_1}f=\mu_0^{-1/2}U_1f$ for $t=\alpha$ and $t_j=p_j$, $j=1,2,3,4$. We have
$$
\align
\Dpinner fg &(\Sigma n_j\delta_j)
\\
& =
\frac1q \sum_{m,n=0}^{q-1} e(\tfrac{nm_3+mm_4}q)
\int_{\Bbb R} \conj{ f(x,n,m) } g(x+\beta n_1, n+m_1, m+m_2) e(\beta n_2x) \,dx
\\
&=
\frac{\sqrt2}{q^2} \sum_{m,n=0}^{q-1} e(\tfrac{nm_3+mm_4}q)
\int_{\Bbb R} \conj{h(x) \varphi(n,m) }
e(\tfrac{t_3(n+m_1)+t_4(m+m_2)}q)\, h(x+t+\beta n_1) 
\\
& \hskip2in \cdot \varphi(n+m_1+t_1, m+m_2+t_2) e(\beta n_2x) \,dx
\\
&=
\frac{\sqrt2}{q^2} H(t+\beta n_1, \beta n_2) \Omega_g
\endalign
$$
where $H(\cdot,\cdot)$ was given in Section 2 and (using the summation in the
proof of Lemma 4.1)
$$
\align
\Omega_g &= 
\sum_{m,n=0}^{q-1} e(\tfrac{nm_3+mm_4}q) e(\tfrac{t_3(n+m_1)+t_4(m+m_2)}q)
\conj{\varphi(n,m) } \varphi(n+m_1+t_1, m+m_2+t_2).
\\
&= e(\tfrac{t_3m_1+t_4m_2}q)
\sum_{m,n=0}^{q-1} e(\tfrac{n(m_3+t_3)+m(m_4+t_4)}q)
\conj{\varphi(n,m) } \varphi(n+m_1+t_1, m+m_2+t_2)
\\
&=
q^2 \, e(\tfrac{t_3m_1+t_4m_2}q)
e(\tfrac1q[a(m_1+t_1)^2+b(m_1+t_1)(m_2+t_2)+\gamma(m_2+t_2)^2]) 
\\
&\ \ \ \ \ \ \ \ \ \ \ \cdot
\delta_q^{(m_3+t_3)+2a(m_1+t_1)+b(m_2+t_2)} 
\delta_q^{(m_4+t_4)+2\gamma(m_2+t_2)+b(m_1+t_1)}
\\
&=
q^2 \, e(\tfrac{t_3m_1+t_4m_2}q)
e(\tfrac1q[a(m_1+t_1)^2+b(m_1+t_1)(m_2+t_2)+\gamma(m_2+t_2)^2])
\\
&\ \ \ \ \ \ \ \ \ \ \ \cdot
\delta_q^{m_3+2am_1+bm_2 + u_3} 
\delta_q^{m_4+2\gamma m_2+bm_1 + u_4}
\endalign
$$
where
$$
u_3 = t_3+2at_1+bt_2, \qquad u_4 = t_4+2\gamma t_2+bt_1.
\tag5.1
$$
Now $\Omega_g$ is non-zero precisely when
$$
m_3+2am_1+bm_2 + u_3 \equiv 0, \qquad m_4+2\gamma m_2+bm_1 + u_4 \equiv 0
$$
which is the diophantine system considered in Section 3, and which by Proposition 3.1
has a unique solution (mod $q$) for $n_3,n_4$.  From (3.4) we obtained
$$
n_3 = c_3 + a_1n_1 + a_2n_2, \qquad n_4 = c_4 + b_1n_1 + b_2n_2,
\tag5.2
$$
where $a_j,b_j,c_j$ are given by (3.5). Thus by (2.4) and (5.2) one gets
$$
\align
\Dpinner fg &= 
\sum_{n_1n_2n_3n_4} \Dpinner fg (\Sigma n_j\delta_j) \Lambda_{n_1n_2n_3n_4}\,
V_4^{n_4} V_3^{n_3} V_2^{n_2} V_1^{n_1}
\\
&= 
\sum_{n_1n_2} C_{n_1n_2} \Lambda_{n_1n_2}'\,
V_4^{c_4+b_1n_1+b_2n_2} V_3^{c_3+a_1n_1+a_2n_2} V_2^{n_2} V_1^{n_1}
\endalign
$$
where, as $m_j$ are given by (2.2), 
$$
\align
C_{n_1n_2} &= \Dpinner fg (\Sigma n_j\delta_j) 
\Big|_{n_3=c_3 + a_1n_1 + a_2n_2,\  n_4=c_4 + b_1n_1 + b_2n_2}
\\
&= \sqrt2 H(t+\beta n_1, \beta n_2) e(\tfrac{t_3m_1+t_4m_2}q) 
e(\tfrac1q[a(m_1+t_1)^2+b(m_1+t_1)(m_2+t_2)+\gamma(m_2+t_2)^2])
\\
&=
e(\tfrac{-t_3t_1-t_4t_2}q) 
e(-\tfrac12(t+\beta n_1)\beta n_2) e^{-\tfrac\pi2[(t+\beta n_1)^2+\beta^2n_2^2]} 
e(\tfrac{t_3(m_1+t_1)+t_4(m_2+t_2)}q) 
\\
&\ \ \ \ \ \ \ \ \ \cdot e(\tfrac1q[a(m_1+t_1)^2+b(m_1+t_1)(m_2+t_2)+\gamma(m_2+t_2)^2])
\endalign
$$
in which the $n_3,n_4$ in $m_1,m_2$ must be evaluated according to (5.2), and by (2.3)
$$
\align
\Lambda_{n_1n_2}' &= \Lambda_{n_1n_2n_3n_4} 
\Big|_{n_3=c_3 + a_1n_1 + a_2n_2,\  n_4=c_4 + b_1n_1 + b_2n_2} 
\\
&= \lambda_{13}^{n_1n_3} \lambda_{24}^{n_2n_4}
\lambda_{12}^{n_1n_2} \lambda_{34}^{n_3n_4}
\lambda_{11}^{n_1(n_1+q)/2}\, \lambda_{22}^{n_2(n_2+q)/2}\,
\lambda_{33}^{n_3(n_3+q)/2}\, \lambda_{44}^{n_4(n_4+q)/2}
\\
&=
\lambda_{13}^{n_1(c_3 + a_1n_1 + a_2n_2)} \lambda_{24}^{n_2(c_4 + b_1n_1 + b_2n_2)}
\lambda_{12}^{n_1n_2} \lambda_{34}^{(c_3 + a_1n_1 + a_2n_2)(c_4 + b_1n_1 + b_2n_2)}
\lambda_{11}^{n_1(n_1+q)/2}
\\
&\ \ \ \ \ \ \ \ \cdot \lambda_{22}^{n_2(n_2+q)/2}\,
\lambda_{33}^{(c_3 + a_1n_1 + a_2n_2)(c_3 + a_1n_1 + a_2n_2+q)/2}\, 
\lambda_{44}^{(c_4 + b_1n_1 + b_2n_2)(c_4 + b_1n_1 + b_2n_2+q)/2}.
\endalign
$$
Let us write
$$
m_1 + t_1 = d_0 + d_1n_1 + d_2n_2, \qquad m_2 + t_2 = e_0 + e_1n_1 + e_2n_2,
$$
where
$$
\aligned
d_0 &= t_1 + p_2c_3 + p_4c_4, \qquad d_1 = -cp_1 + p_2a_1 + p_4b_1, \qquad
d_2 = cp_3 + p_2a_2 + p_4b_2
\\
e_0 &= t_2 - p_1c_3 - p_3c_4, \qquad e_1 = -cp_2 - p_1a_1 - p_3b_1, \qquad
e_2 = cp_4 - p_1a_2 - p_3b_2
\endaligned
\tag5.3
$$
Thus
$$
\align
C_{n_1n_2} &= e(\tfrac{-t_3t_1-t_4t_2}q) 
e(-\tfrac12(t+\beta n_1)\beta n_2) e^{-\tfrac\pi2[(t+\beta n_1)^2+\beta^2n_2^2]} 
e(\tfrac{t_3(d_0 + d_1n_1 + d_2n_2)+t_4(e_0 + e_1n_1 + e_2n_2)}q)
\\
&
\cdot e(\tfrac1q[a(d_0 + d_1n_1 + d_2n_2)^2+b(d_0 + d_1n_1 + d_2n_2)(e_0 + e_1n_1 + e_2n_2)
+\gamma(e_0 + e_1n_1 + e_2n_2)^2])
\endalign
$$
Now it can be checked that
$$
\multline
V_4^{c_4+b_1n_1+b_2n_2} V_3^{c_3+a_1n_1+a_2n_2} V_2^{n_2} V_1^{n_1}
\\
=
e(-\tfrac pq c_3(b_1n_1+b_2n_2)) e(-\tfrac pq b_1a_2n_1n_2) 
e(-\tfrac pq a_2b_2 \tfrac{n_2(n_2-1)}2) e(-\tfrac pq a_1b_1 \tfrac{n_1(n_1-1)}2) 
V_4^{c_4}V_3^{c_3} X_2^{n_2} X_1^{n_1}
\endmultline
$$
where $X_1 = V_4^{b_1}V_3^{a_1} V_1,\ X_2 = V_4^{b_2}V_3^{a_2} V_2$ are unitaries that
generate some rotation algebra. Letting
$$
R = 2c_3(b_1n_1+b_2n_2) + 2b_1a_2n_1n_2 + a_2b_2n_2(n_2-1) + a_1b_1n_1(n_1-1)
$$
one has
$$
V_4^{c_4+b_1n_1+b_2n_2} V_3^{c_3+a_1n_1+a_2n_2} V_2^{n_2} V_1^{n_1}
= e(-\tfrac {pR}{2q}) V_4^{c_4}V_3^{c_3} X_2^{n_2} X_1^{n_1}.
$$
Using the following data (recall $\lambda_{jk}=\frak h(\delta_j,\delta_k)$ and 
$C=p_1p_4-p_2p_4$)
$$
\align
\lambda_{13} &= e(\tfrac cq(p_1p_4-p_2p_3) = e(\tfrac{cC}q)
\\
\lambda_{24} &= e(-\tfrac{cC}q) = \lambda_{13}^{-1}
\\
\lambda_{12} &= e(\beta^2) e(\tfrac{c^2(p_1^2+p_2^2)}q) = e(\beta^2) \lambda_{34}^{c^2}
\\
\lambda_{34} &= e(\tfrac{p_1^2+p_2^2}q)
\\
\lambda_{11} &= e(\tfrac{c^2(p_1p_3+p_2p_4)}q) = \lambda_{44}^{c^2}
\\
\lambda_{22} &= \lambda_{44}^{-c^2}
\\
\lambda_{33} &= \lambda_{44}^{-1}
\\
\lambda_{44} &= e(\tfrac{p_1p_3+p_2p_4}q)
\endalign
$$
one gets
$$
\align
&\Lambda_{n_1n_2}' 
\\
&= e(\beta^2n_1n_2)
\lambda_{13}^{n_1(c_3 + a_1n_1 + a_2n_2) - n_2(c_4 + b_1n_1 + b_2n_2)}
\lambda_{34}^{c^2n_1n_2 + (c_3 + a_1n_1 + a_2n_2)(c_4 + b_1n_1 + b_2n_2)}
\\
&\cdot \lambda_{44}^{c^2n_1(n_1+q)/2-c^2n_2(n_2+q)/2 - 
(c_3 + a_1n_1 + a_2n_2)(c_3 + a_1n_1 + a_2n_2+q)/2 + 
(c_4 + b_1n_1 + b_2n_2)(c_4 + b_1n_1 + b_2n_2+q)/2}
\endalign
$$
which we write as
$$
\Lambda_{n_1n_2}' = e(\beta^2n_1n_2) e(\tfrac L{2q})
$$
where $L$ is the integer
$$
\align
L &= 2cC \bigr[n_1(c_3 + a_1n_1 + a_2n_2) - n_2(c_4 + b_1n_1 + b_2n_2)\bigr]
\\
&\ + 
2(p_1^2+p_2^2)\bigl[ c^2n_1n_2 +  (c_3 + a_1n_1 + a_2n_2)(c_4 + b_1n_1 + b_2n_2)\bigr]
\\
&\ + 
(p_1p_3+p_2p_4) \Bigl[c^2n_1(n_1+q) - c^2n_2(n_2+q) - 
(c_3 + a_1n_1 + a_2n_2)(c_3 + a_1n_1 + a_2n_2+q)
\\
&\ \ \ \ \ \ +
(c_4 + b_1n_1 + b_2n_2)(c_4 + b_1n_1 + b_2n_2+q) \Bigr].
\endalign
$$
Let us also write
$$
C_{n_1n_2} = e^{-\tfrac\pi2[(t+\beta n_1)^2+\beta^2n_2^2]} 
e(-\tfrac12\beta^2 n_1n_2) e(-\tfrac12 t\beta n_2) \, e(\tfrac Gq)
$$
where $G$ is the integer
$$
\align
G &= -t_3t_1-t_4t_2 + t_3(d_0 + d_1n_1 + d_2n_2)+t_4(e_0 + e_1n_1 + e_2n_2)
\\
&\ \ \ 
+ a(d_0 + d_1n_1 + d_2n_2)^2+b(d_0 + d_1n_1 + d_2n_2)(e_0 + e_1n_1 + e_2n_2)
+\gamma(e_0 + e_1n_1 + e_2n_2)^2.
\endalign
$$
Therefore,
$$
\Dpinner fg = V_4^{c_4}  V_3^{c_3}
\sum_{n_1,n_2}  e^{-\tfrac\pi2[(t+\beta n_1)^2+\beta^2n_2^2]} 
e(\tfrac12\beta^2 n_1n_2) e(-\tfrac12 t\beta n_2)
e(\tfrac S{2q}) X_2^{n_2} X_1^{n_1}
$$
where $S = 2G+L-pR$ is the integer
$$
\align
S &= -2t_3t_1-2t_4t_2 + 2t_3(d_0 + d_1n_1 + d_2n_2)+2t_4(e_0 + e_1n_1 + e_2n_2)
\\
&\ \ \ 
+ 2a(d_0 + d_1n_1 + d_2n_2)^2 + 2b(d_0 + d_1n_1 + d_2n_2)(e_0 + e_1n_1 + e_2n_2)
+ 2\gamma(e_0 + e_1n_1 + e_2n_2)^2
\\
&\ \ \ + 2cC \bigr[n_1(c_3 + a_1n_1 + a_2n_2) - n_2(c_4 + b_1n_1 + b_2n_2)\bigr]
\\
&\ \ \ + 
2(p_1^2+p_2^2)\bigl[ c^2n_1n_2 +  (c_3 + a_1n_1 + a_2n_2)(c_4 + b_1n_1 + b_2n_2)\bigr]
\\
&\ \ \ + 
(p_1p_3+p_2p_4) \Bigl[c^2n_1(n_1+q) - c^2n_2(n_2+q) - 
(c_3 + a_1n_1 + a_2n_2)(c_3 + a_1n_1 + a_2n_2+q)\Bigr.
\\
&\ \ \ \ \ \ \ \ \ \ \ \ \ +
\Bigl.(c_4 + b_1n_1 + b_2n_2)(c_4 + b_1n_1 + b_2n_2+q) \Bigr]
\\
&\ \ \  -2pc_3(b_1n_1+b_2n_2) - 2pb_1a_2n_1n_2 - pa_2b_2n_2(n_2-1) - pa_1b_1n_1(n_1-1).
\endalign
$$
Accumulating coefficients one writes this in the form
$$
S= a'n_1^2 + b'n_2^2 + 2c'n_1n_2 + d'n_1 + e'n_2 + K
$$
where 
$$
a' = 2ad_1^2 + 2bd_1e_1 + 2\gamma e_1^2 + 2cCa_1 + 2(p_1^2+p_2^2)a_1b_1 +
(p_1p_3+p_2p_4)[c^2 - a_1^2 + b_1^2] - pa_1b_1
$$
$$
b' = 2ad_2^2 + 2bd_2e_2 + 2\gamma e_2^2 - 2cCb_2 + 2(p_1^2+p_2^2)a_2b_2 +
(p_1p_3+p_2p_4)[-c^2 - a_2^2 + b_2^2] - pa_2b_2
$$
$$
\align
c' &= 2ad_1d_2 + b(d_1e_2+d_2e_1) + 2\gamma e_1e_2 + cC(a_2 - b_1) +
(p_1^2+p_2^2)[c^2 + a_1b_2 + a_2b_1] 
\\ &\ \ \ \ \ \ \ \ \ \ +(p_1p_3+p_2p_4)(b_1b_2-a_1a_2) - pa_2b_1
\\
d' &= 2t_3d_1 + 2t_4e_1 + 4ad_0d_1 + 2b(d_0e_1+d_1e_0) + 4\gamma e_0e_1 + 
2cCc_3 + 2(p_1^2+p_2^2)[ c_3b_1 + c_4a_1 ]
\\
&\ \ \ \ \ \ \ \ 
+ (p_1p_3+p_2p_4)[qc^2-(c_3a_1 + a_1(c_3+q)) + (c_4b_1 + b_1(c_4+q))] - 2pc_3b_1 + pa_1b_1
\\
e' &= 2t_3d_2 + 2t_4e_2 + 4ad_0d_2 + 2b(d_0e_2+d_2e_0) + 4\gamma e_0e_2 -
2cCc_4 + 2(p_1^2+p_2^2)[ c_3b_2 + c_4a_2  ]
\\
&\ \ \ \ \ \ \ \ 
+ (p_1p_3+p_2p_4)[-qc^2-(c_3a_2 + a_2(c_3+q)) + (c_4b_2 + b_2(c_4+q))] - 2pc_3b_2 + pa_2b_2
\endalign
$$
and
$$
\align
K &= -2t_3t_1-2t_4t_2 + 2t_3d_0 + 2t_4e_0 + 2ad_0^2 + 2bd_0e_0 + 2\gamma e_0^2 
+ 2(p_1^2+p_2^2)c_3c_4
\\
&\ \ \ \ \ \ \ \ + (p_1p_3+p_2p_4)[ -c_3(c_3+q) + c_4(c_4+q) ].
\endalign
$$
Let us write
$$
d' = 2d_0' + d_1', \qquad e' = 2e_0' + e_1'
$$
where $d_0',e_0'$ consist of the terms that depend on $t_j$, and $d_1', e_1'$ consist
of terms that do not.  Thus
$$
\align
d_0' &= t_3d_1 + t_4e_1 + 2ad_0d_1 + b(d_0e_1+d_1e_0) + 2\gamma e_0e_1 + 
cCc_3 + (p_1^2+p_2^2)[ c_3b_1 + c_4a_1 ]
\\
&\ \ \ \ \ \ \ \ 
+ (p_1p_3+p_2p_4)[b_1c_4-a_1c_3] - pc_3b_1,
\\
e_0' &= t_3d_2 + t_4e_2 + 2ad_0d_2 + b(d_0e_2+d_2e_0) + 2\gamma e_0e_2 -
cCc_4 + (p_1^2+p_2^2)[ c_3b_2 + c_4a_2  ]
\\
&\ \ \ \ \ \ \ \ 
+ (p_1p_3+p_2p_4)[b_2c_4-c_3a_2] - pc_3b_2,
\endalign
$$
and
$$
\align
d_1' &= (p_1p_3+p_2p_4)[qc^2-a_1q + b_1q] + pa_1b_1, 
\\
e_1' &= (p_1p_3+p_2p_4)[-qc^2-a_2q + b_2q] + pa_2b_2.
\endalign
$$
In the following lemma it will be shown that $c' \equiv 0 \mod q$, hence we have
$$
\multline
\Dpinner fg = e(\tfrac K{2q})  V_4^{c_4}  V_3^{c_3}
\sum_{n_1,n_2}  e^{-\tfrac\pi2[(t+\beta n_1)^2+\beta^2n_2^2]} 
e(\tfrac12\beta^2 n_1n_2) e(-\tfrac12 t\beta n_2) 
\\
\cdot e(\tfrac{a'}{2q}n_1^2 + \tfrac{b'}{2q}n_2^2 + \tfrac{d'}{2q}n_1 + \tfrac{e'}{2q}n_2) X_2^{n_2} X_1^{n_1}
\endmultline
\tag5.4
$$

\medpagebreak

\proclaim{Lemma 5.1} One has
$a' \equiv b' \equiv c' \equiv d_0' \equiv e_0' \equiv 0 \mod q$.
\endproclaim

\demo{Proof}
Let $r = \tfrac{a'}{2q}, \ s= \tfrac{b'}{2q}$. We shall first prove that $r,s$ must be
half integers, so it will follow that $a' \equiv b' \equiv 0$.  Specializing (5.4) to the
case $g=f$, so that $t=t_j=0$ and $K=c_3=c_4=0$ (and noting that $a',b'$ are the same for
$g=f$ and $g=\pi_{\varepsilon_1}f$), one gets
$$
\Dpinner ff = \sum_{m,n}  e^{-\tfrac\pi2\beta^2(m^2+n^2)} 
e(\tfrac12\beta^2 mn) e(rm^2 + sn^2) Y_2^n Y_1^m
$$
where $Y_1 = e(\tfrac{d'}{2q})X_1$ and $Y_2 = e(\tfrac{e'}{2q})X_2$ are 
unitaries that generate some rotation algebra (each of whose spectrum is the unit circle).
Write 
$$
\Dpinner ff = \sum_m e(rm^2) e^{-\tfrac\pi2 \beta^2 m^2} \, \rho_m \ Y_1^m
$$
where 
$$
\rho_m = \sum_n e(sn^2) e^{-\tfrac\pi2 \beta^2 n^2} e(\tfrac12\beta^2mn) Y_2^n.
$$
Viewing $\rho_m$ as a period 1 function of the real variable $t$, it can be written as
$$
\rho_m(t) \ = \ 
\sum_n e(sn^2) e^{-\tfrac\pi2 \beta^2 n^2} \ e(\tfrac12\beta^2mn) e(nt)
\ = \ \vartheta_3(\pi t + \tfrac\pi2\beta^2 m,\, \tfrac i2\beta^2+2s).
$$
Since $\Dpinner ff$ is a positive element, its conditional expectation $\rho_0$ is also
positive. (Here, the conditional expectation is the map taking a generic sum
$\sum_m f_m Y_1^m$ to $f_0$, where $f_m$ are functions of $Y_2$.)  
Since $\rho_0(t)$ must in particular be a real function, this forces $s$ to be a half integer.
To see this, one simply sets $\rho_0(t)$ equal to its conjugate:
$$
\sum_n e(sn^2) e^{-\tfrac\pi2 \beta^2 n^2} \ e(nt)
= \sum_n e(-sn^2) e^{-\tfrac\pi2 \beta^2 n^2} \ e(-nt)
= \sum_n e(-sn^2) e^{-\tfrac\pi2 \beta^2 n^2} \ e(nt).
$$
Since this holds for all real $t$, one must have $e(sn^2) = e(-sn^2)$ for all $n$.
Putting $n=1$ one gets the result.
By symmetry, one similarly deduces that $r$ is a half integer.

Since $d_0'=e_0'=0$ trivially in the case $g=f$, we focus on the case $g=\pi_{\varepsilon_1}f$
and insert $t_j=p_j$ in the above expressions for $d_0',e_0'$. Let us first prove that  
$d_0' \equiv 0 \mod q$ where
$$
\align
d_0' &= p_3d_1 + p_4e_1 + 2ad_0d_1 + b(d_0e_1+d_1e_0) + 2\gamma e_0e_1 + 
c(p_1p_4-p_2p_3)c_3
\\
&\ \ \ \ \ \ \ \ 
+ (p_1^2+p_2^2)(c_3b_1 + c_4a_1) + (p_1p_3+p_2p_4)[b_1c_4-a_1c_3] - pc_3b_1,
\endalign
$$
where $d_0 = p_1 + p_2c_3 + p_4c_4,\ e_0 = p_2 - p_1c_3 - p_3c_4$, and in $c_3,c_4$ one
puts $u_3=s_3,\ u_4=s_4$ (see (3.5) and (5.1)).
From the definitions (3.5), we have $a_1=-cc_3,\ b_1=-cc_4$.  Also, it is 
straighforward to check that $d_1=-cd_0,\ e_1=-ce_0$.  Putting $n_1=n_2=0$ in (3.4) and 
inserting into (3.2) gives us the congruences (mod $q$)
$$
r_4c_3 + r_2c_4 \equiv -s_3, \qquad r_3c_3 + r_1c_4 \equiv -s_4. 
\tag5.5
$$
Using the congruence $pc\equiv1$, we rewrite $d_0'$ as
$$
\align
d_0' &\equiv -cd_0(p_3+2ad_0+be_0) - ce_0(p_4+2\gamma e_0+bd_0) 
\\
&\ \ + c(p_1p_4-p_2p_3)c_3 - 2c(p_1^2+p_2^2)c_3c_4
- c(p_1p_3+p_2p_4)(c_4^2-c_3^2) + pcc_3c_4,
\endalign
$$
Using the congruences (5.5) one verifies that
$$
p_3+2ad_0+be_0 \equiv p_4c_3 - p_2c_4, \qquad  p_4+2\gamma e_0+bd_0 \equiv p_1c_4 - p_3c_3.
\tag5.6
$$
Thus
$$
\align
d_0' &\equiv -cd_0(p_4c_3 - p_2c_4) - ce_0(p_1c_4 - p_3c_3) 
\\
&\ \ + c(p_1p_4-p_2p_3)c_3 - 2c(p_1^2+p_2^2)c_3c_4
- c(p_1p_3+p_2p_4)(c_4^2-c_3^2) + pcc_3c_4
\\
&= -c(p_1+p_2c_3+p_4c_4)(p_4c_3 - p_2c_4) - c(p_2-p_1c_3-p_3c_4)(p_1c_4 - p_3c_3) 
\\
&\ \ + c(p_1p_4-p_2p_3)c_3 - 2c(p_1^2+p_2^2)c_3c_4
- c(p_1p_3+p_2p_4)(c_4^2-c_3^2) + pcc_3c_4
\endalign
$$
factoring out $c$ and expanding this last expression one checks that all of its terms 
formally cancel out and give zero (after inserting $p = p_1^2+p_2^2+p_3^2+p_4^2$).

\medpagebreak

In a similar manner we show that $e_0' \equiv 0$, where
$$
\align
e_0' &= p_3d_2 + p_4e_2 + 2ad_0d_2 + b(d_0e_2+d_2e_0) + 2\gamma e_0e_2 -
c(p_1p_4-p_2p_3)c_4 
\\
&\ \ \ \ \ \ \ \ 
+ (p_1^2+p_2^2)(c_3b_2 + c_4a_2) + (p_1p_3+p_2p_4)[b_2c_4-c_3a_2] - pc_3b_2.
\endalign
$$
Since $c_3 \equiv -pa_1,\ c_4 \equiv -pb_1$, using (5.6), and $pc\equiv1$, one has
$$
\align
e_0' &\equiv d_2(p_3+2ad_0+be_0) + e_2(p_4+2\gamma e_0+bd_0) 
+ (p_1p_4-p_2p_3)b_1  - p(p_1^2+p_2^2)(a_1b_2 + a_2b_1)
\\
&\ \ \ \ \ \ \ \ 
- p(p_1p_3+p_2p_4)(b_1b_2-a_1a_2) + p^2a_1b_2
\\
&\equiv -p(cp_3+p_2a_2+p_4b_2)(p_4a_1-p_2b_1) - p(cp_4-p_1a_2-p_3b_2)(p_1b_1-p_3a_1) 
\\
&\ \ \ \ \ \ \ \ + (p_1p_4-p_2p_3)b_1  - p(p_1^2+p_2^2)(a_1b_2 + a_2b_1)
- p(p_1p_3+p_2p_4)(b_1b_2-a_1a_2) + p^2a_1b_2
\\
&\equiv -p_3(p_4a_1-p_2b_1) - p(p_2a_2+p_4b_2)(p_4a_1-p_2b_1) 
- p_4(p_1b_1-p_3a_1)  
\\
&\ \ \ \ \ \ \ \ + p(p_1a_2+p_3b_2)(p_1b_1-p_3a_1) 
+ (p_1p_4-p_2p_3)b_1  - p(p_1^2+p_2^2)(a_1b_2 + a_2b_1)
\\
&\ \ \ \ \ \ \ \ - p(p_1p_3+p_2p_4)(b_1b_2-a_1a_2) + p^2a_1b_2
\endalign
$$
or
$$
\align
e_0' &\equiv
- p(p_2a_2+p_4b_2)(p_4a_1-p_2b_1) + p(p_1a_2+p_3b_2)(p_1b_1-p_3a_1) 
- p(p_1^2+p_2^2)(a_1b_2 + a_2b_1)
\\
&\ \ \ \ \ \ \ \ - p(p_1p_3+p_2p_4)(b_1b_2-a_1a_2) + p^2a_1b_2
\endalign
$$
factoring out $p$ and expanding this last expression one checks that all of its terms 
formally cancel out and give zero.

\medpagebreak

Now we prove that $c' \equiv 0 \mod q$.  Let us write
$$
\align
c' &= Q + c(p_1p_4-p_2p_3)(a_2 - b_1) +
(p_1^2+p_2^2)[c^2 + a_1b_2 + a_2b_1] 
\\ &\ \ \ \ \ \ \ \ \ \ +(p_1p_3+p_2p_4)(b_1b_2-a_1a_2) - pa_2b_1
\endalign
$$
where, using (5.6) and freely inserting $pc\equiv1$ when needed,
$$
\align
Q &:= 2ad_1d_2 + b(d_1e_2+d_2e_1) + 2\gamma e_1e_2
\\
&= d_2(2ad_1+be_1) + e_2(2\gamma e_1 + bd_1)
\\
&= -cd_2(2ad_0+be_0) - ce_2(2\gamma e_0 + bd_0)
\\
&\equiv -cd_2(-p_3+p_4c_3-p_2c_4) - ce_2(-p_4+p_1c_4-p_3c_3)
\\
&= d_2(cp_3+p_4a_1-p_2b_1) + e_2(cp_4+p_1b_1-p_3a_1)
\\
\intertext{inserting the values for $d_2,e_2$ gives}
&= (cp_3+p_2a_2+p_4b_2)cp_3 + (cp_4-p_1a_2-p_3b_2)cp_4 
\\
&\ \ \ \ \ \ \ + (cp_3+p_2a_2+p_4b_2)(p_4a_1-p_2b_1) + (cp_4-p_1a_2-p_3b_2)(p_1b_1-p_3a_1)
\\
&= (cp_3+p_2a_2+p_4b_2)cp_3 + (cp_4-p_1a_2-p_3b_2)cp_4 
+ cp_3(p_4a_1-p_2b_1) + cp_4(p_1b_1-p_3a_1)
\\
&\ \ \ \ \ \ \ + (p_2a_2+p_4b_2)(p_4a_1-p_2b_1) + (-p_1a_2-p_3b_2)(p_1b_1-p_3a_1)
\\
&= c^2p_3^2 + c^2p_4^2 + cp_2p_3a_2 - cp_1p_4a_2 - cp_2p_3b_1 + cp_1p_4b_1
\\
&\ \ \ \ \ \ \ + (p_2a_2+p_4b_2)(p_4a_1-p_2b_1) - (p_1a_2+p_3b_2)(p_1b_1-p_3a_1)
\endalign
$$
after expanding and making the several cancellations we get
$$
c' \equiv c^2 p + p(a_1b_2 - a_2b_1).
$$
Substituting the values for $a_j,b_j$ (see (3.5)) gives
$$
\align
c' &\equiv c^2 p + c^2p\,(\Delta')^2
\left[(r_1s_3-r_2s_4)(r_4s_2-r_3s_1) - (r_1s_1-r_2s_2)(r_4s_4-r_3s_3) \right] 
\\
&= c^2 p - c^2p\,(\Delta')^2(r_1r_4-r_2r_3)(s_1s_4-s_2s_3)
\\
&= c^2 p - c^2p\,(\Delta')^2 \Delta^2
\\
&\equiv 0
\endalign
$$
since (from Section 3) $\Delta = r_1r_4-r_2r_3 = s_1s_4-s_2s_3$ and $\Delta'\Delta\equiv 1$.
\qed
\enddemo

\newpage

Therefore, by Lemma 5.1, and writing $a'=qa'', b'=qb''$ for some integers $a'',b''$ 
(which are independent of $t_j$, as $a',b'$ are independent), one gets from (5.4) 
(for $g=f$ and $g=\pi_{\varepsilon_1}f$)
$$
\multline
\Dpinner fg = e(\tfrac K{2q})  V_4^{c_4}  V_3^{c_3}
\sum_{n_1,n_2}  e^{-\tfrac\pi2[(t+\beta n_1)^2+\beta^2n_2^2]} 
e(\tfrac12\beta^2 n_1n_2) e(-\tfrac12 t\beta n_2) 
\\
\cdot e(\tfrac{a''}2n_1^2 + \tfrac{b''}2n_2^2 + \tfrac{d_1'}{2q}n_1 + \tfrac{e_1'}{2q}n_2) X_2^{n_2} X_1^{n_1}
\endmultline
$$
Since $e(\tfrac{a''}2) = \pm1 = e(\tfrac{b''}2)$ one can write 
$e(\tfrac{a''}2n_1^2) = e(\tfrac{a''}2n_1)$ so that
$$
\Dpinner fg = e(\tfrac K{2q})  V_4^{c_4}  V_3^{c_3}
\sum_{m,n}  e^{-\tfrac\pi2[(t+\beta m)^2+\beta^2n^2]} 
e(\tfrac12\beta^2 mn) e(-\tfrac12 t\beta n) W_2^n W_1^m
$$
where
$$
\align
W_1 &= e(\tfrac{a''}2+\tfrac{d_1'}{2q}) X_1 
= e(\tfrac{a''}2+\tfrac{d_1'}{2q}) V_4^{b_1}V_3^{a_1} V_1,
\\
W_2 &= e(\tfrac{b''}2+\tfrac{e_1'}{2q}) X_2 
= e(\tfrac{b''}2+\tfrac{e_1'}{2q}) V_4^{b_2}V_3^{a_2} V_2,
\endalign 
$$
are unitaries that do not depend on $t_j$, but which generate some rotation 
algebra (each of whose spectrum is the unit circle).  Specializing this to the 
cases $g=f$ and $g=\pi_{\varepsilon_1}f = \mu_0^{-1/2}U_1f$ one gets
$$
\Dpinner ff = \sum_{m,n}  e^{-\tfrac\pi2\beta^2(m^2+n^2)} 
e(\tfrac12\beta^2 mn)  W_2^n W_1^m
\tag5.7
$$
and (recalling that $\beta =\tfrac1{q\alpha}$)
$$
\Dpinner f{U_1f} = \mu_0^{1/2} e(\tfrac K{2q})  V_4^{c_4}  V_3^{c_3}
\sum_{m,n}  e^{-\tfrac\pi2[(\alpha+\beta m)^2+\beta^2n^2]} 
e(\tfrac12\beta^2 mn) e(-\tfrac n{2q}) W_2^n W_1^m
\tag5.8
$$
where in the latter case one has $t_j=p_j$ and, by (3.5),
$$
c_3 = \Delta'(-r_1s_3+r_2s_4), \qquad c_4 = \Delta'(-r_4s_4+r_3s_3),
$$
since, from (5.1), $u_3  =  p_3+2ap_1+bp_2  =  s_3$ and 
$u_4  =  p_4+2\gamma p_2+bp_1  =  s_4$.

\newpage

\subhead \S6. INVERTIBILTY OF $\Dpinner ff$ \endsubhead

From (5.7) we obtained
$$
\Dpinner ff = \sum_{m,n}  e^{-\tfrac\pi2\beta^2(m^2+n^2)} 
e(\tfrac12\beta^2 mn)  W_2^n W_1^m
$$
and our objective in this section is to show that (because $\beta^2>1$) this positive
element is always invertible in the rotation algebra generated by $W_1,W_2$-----in fact
one can easily check that $W_1W_2 = e(\beta^2)W_2W_1$.  See Theorem 6.4 below.  As in 
the proof of Lemma 5.1 we write
$$
\Dpinner ff = \sum_m e^{-\tfrac\pi2 \beta^2 m^2} \, \rho_m \ W_1^m
$$
where 
$$
\rho_m = \sum_n e^{-\tfrac\pi2 \beta^2 n^2} e(\tfrac12\beta^2mn) W_2^n.
$$
Viewing $\rho_m$ as a period 1 function on $\Bbb R$ (as $W_2$ has spectrum the 
unit circle), it can be written as
$$
\rho_m(t) \ = \ 
\sum_n e^{-\tfrac\pi2 \beta^2 n^2} \ e(\tfrac12\beta^2mn) e(nt)
\ = \ \vartheta_3(\pi t + \tfrac\pi2\beta^2 m,\, \tfrac i2\beta^2).
$$
We have $\rho_0(t) = \vartheta_3(\pi t,\, \tfrac i2\beta^2)$ and it is
well-known from Theta function theory that for all real $t$ 
$$
\vartheta_3(0,\, \tfrac i2\beta^2) 
\ge \vartheta_3(\pi t,\, \tfrac i2\beta^2) 
\ge \vartheta_3(\tfrac\pi2,\, \tfrac i2\beta^2) > 0
$$
so that $\rho_0>0$. Note that in fact
$\rho_m(t) = \rho_0(t+\tfrac12\beta^2 m)$ 
so that also each $\rho_m$ is also strictly positive.  In particular, we have the norms
$$
\|\rho_m\| = \|\rho_0\| = \vartheta_3(0,\, \tfrac i2\beta^2), \qquad
\|\rho_m^{-1}\| = \|\rho_0^{-1}\| = \frac1{\vartheta_3(\tfrac\pi2,\, \tfrac i2\beta^2)}
\tag6.1
$$
for all $m$.  

\bigpagebreak

\proclaim{Lemma 6.1}
If $f(t) = \vartheta_3(\pi t + 2\pi bm,ib)$ and $g(t) = \vartheta_3(\pi t,ib)$, where
$b>\tfrac12$ and $m\in\Bbb Z$, then
$$
\|f-g\| \ \le \ 8\pi |m| (b-\tfrac12) \sum_{k=1}^\infty k e^{-\pi bk^2}.
$$
\endproclaim
\demo{Proof} By definition of $\vartheta_3$ we have
$$
f(t)-g(t) = \sum_k e^{-\pi bk^2} e^{2\pi i tk} [e(2bmk)-1].
$$
Since $e(2bmk)=e(2(b-\tfrac12)mk)$ and $|e(x)-1| \le 2\pi|x|$ for all $x$, one gets
$$
\align
|f(t)-g(t)| &\le \sum_k e^{-\pi bk^2} |e(2(b-\tfrac12)mk)-1|
\le 4\pi(b-\tfrac12)|m| \sum_k |k| e^{-\pi bk^2}
\endalign
$$
hence the result. \qed
\enddemo
It now follows from Lemma 6.1 that for each integer $m$
$$
\|\rho_{2m} - \rho_0\| \le 
4\pi|m|(\beta^2-1) \sum_{k=1}^\infty k e^{-\tfrac\pi2 \beta^2k^2}.
\tag6.2
$$
For any $m$ we have
$$
\|\rho_{2m+1} - \rho_0\| \le
\vartheta_3(0,\tfrac i2\beta^2) - \vartheta_3(\tfrac\pi2,\tfrac i2\beta^2)
= 2\vartheta_2(0,2i\beta^2).
\tag6.3
$$
We will use these norm bounds shortly.  Consider the operator
$$
R := \sum_m e^{-\tfrac\pi2 \beta^2 m^2} \, W_1^m.
$$
Since, as a period 1 function of the real variable $t$, $R(t)$ is a Theta function 
similar to $\rho_0$,  $R$ is positive and invertible (in the C*-algebra generated by $W_1$).
One has
$$
\| \Dpinner ff - \rho_0 R \| = 
\left\| \sum_m e^{-\tfrac\pi2 \beta^2 m^2} (\rho_m-\rho_0) \, W_1^m \right\|
\le 
\sum_m e^{-\tfrac\pi2 \beta^2 m^2} \|\rho_m-\rho_0\|.
$$
Now divide this sum according to the parity of $m$ using the norm bounds above for even
$m$ and odd $m$.  One gets
$$
\sum_m e^{-\tfrac\pi2 \beta^2 m^2} \|\rho_m-\rho_0\| \ = \
\sum_m e^{-2\pi\beta^2 m^2} \|\rho_{2m}-\rho_0\| +
\sum_m e^{-2\pi \beta^2 (m+\tfrac12)^2} \|\rho_{2m+1}-\rho_0\|.
$$
From (6.2)
$$
\sum_m e^{-2\pi\beta^2 m^2} \|\rho_{2m}-\rho_0\|
\le 
8\pi (\beta^2-1) \Psi(\tfrac12\beta^2) \Psi(2\beta^2) 
$$
where we have written
$$
\Psi(x) = \sum_{k=1}^\infty k e^{-\pi xk^2}
$$
for $x>0$.  Therefore, using (6.3) we obtain
$$
\sum_m e^{-\tfrac\pi2 \beta^2 m^2} \|\rho_m-\rho_0\| \ \le \ 
8\pi (\beta^2-1) \Psi(\tfrac12\beta^2) \Psi(2\beta^2) +
2\vartheta_2(0,2i\beta^2)^2.
$$
(Note:~there is a square on the $\vartheta_2$ here.) Now since 
$\|R^{-1}\| = \|\rho_0^{-1}\| = \vartheta_3(\tfrac\pi2,\tfrac i2\beta^2)^{-1}$ by (6.1),
one obtains
$$
\| R^{-1} \rho_0^{-1} \Dpinner ff  - 1 \| 
\le 
\frac{8\pi (\beta^2-1) \Psi(\tfrac12\beta^2) \Psi(2\beta^2) +
2\vartheta_2(0,2i\beta^2)^2}
{\vartheta_3(\tfrac\pi2,\tfrac i2\beta^2)^2}.
\tag6.4
$$
That the quantity on the right side is less than 1 will now follow from the following 
lemma and the fact that $\beta>1$.

\bigpagebreak

\proclaim{Lemma 6.2} For $x>1$, one has
$$
\frac{8\pi (x-1) \Psi(\tfrac12x) \Psi(2x) + 2\vartheta_2(0,2ix)^2}
{\vartheta_3(\tfrac\pi2,\tfrac i2 x)^2} 
\ < \ 
\frac{2 \vartheta_2(0,2i)^2} {\vartheta_3(\tfrac\pi2,\tfrac i2)^2} \ = \ 1.
$$
\endproclaim
\demo{Proof} 
First note that since $\Psi(x) = e^{-\pi x} \sum_{k\ge1} k e^{-\pi x(k^2-1)}$ one has
$$
\Psi(\tfrac12x) = e^{-\pi x/2} \sum_{k=1}^\infty k e^{-\tfrac12\pi x(k^2-1)}  
\le e^{-\pi x/2} \sum_{k=1}^\infty k e^{-\tfrac12\pi (k^2-1)}  
< 1.01798 e^{-\pi x/2}
$$
for $x>1$.  Similarly,
$$
\Psi(2x) = e^{-2\pi x} \sum_{k=1}^\infty k e^{-2\pi x(k^2-1)}
< 1.000000014 e^{-2\pi x}
$$
(seven zeros after the decimal).  Hence
$$
8\pi (x-1) \Psi(\tfrac12x) \Psi(2x) \le K(x-1)e^{-5\pi x/2}
$$
where $K=8\pi(1.018)$. Therefore it suffice to show that for all $x>1$
$$
\Cal E(x) := \frac{ K(x-1)e^{-5\pi x/2} + 2\vartheta_2(0,2ix)^2}
{\vartheta_3(\tfrac\pi2,\tfrac i2 x)^2} \ < \ 1. 
\tag6.5
$$
The function $(x-1)e^{-5\pi x/2}$ is decreasing for $x> x_0=1+\tfrac2{5\pi} = 1.1273...$.
Since also $\vartheta_2(0,2ix)^2$ is decreasing, while 
$\vartheta_3(\tfrac\pi2,\tfrac i2 x)$ is increasing, it follows that for $x>x_0$ the 
function $\Cal E(x)$ is decreasing and hence $\Cal E(x) \le \Cal E(x_0) < 0.532<1$.

So now we can restrict our attention to $1<x\le1.128$ and show that (6.5) also holds 
for such values.  To do this we will first show that
$$
h(x) := K(x-1) e^{-5\pi x/2} \ < \ 
\vartheta_3(0,2ix)-(1+\sqrt2)\vartheta_2(0,2ix)  =: g(x)
\tag6.6
$$
for $1<x\le1.128$.  This follows by establishing three steps.  First, we use the fact 
that $h$ is underneath its tangent line at $x=1$.  Second, that this tangent line 
is underneath the secant line of $g$ over the interval $[1,1.128]$.  Third, that 
$g$ is above this secant line on $[1,1.128]$.  
The first step is easy to check.
The second step follows since $h'(1) = 0.00993...$ is less 
than the slope of the secant line, which is $[g(1.128)-g(1)]/0.128 = 1.412...$, 
(both lines going through the point $(1,0)$).  The fact that $g(1)=0$, which is needed 
here, is proved in Lemma 6.3 below.  
The third step follows from $g''(x)<0$ for all $x\ge1$, so that $g$ lies above its 
secant line.  To see this, using the definition of Theta functions one checks that
$$
g''(x) \ = \ 8\pi^2 \sum_{k=0}^\infty 
\left( (k+1)^4 e^{-2\pi x(k+1)^2} - (1+\sqrt2)(k+\tfrac12)^4 e^{-2\pi x(k+\tfrac12)^2}
\right).
$$
To see that this is negative one notes that each term of the sum is negative for $x\ge1$ 
and $k\ge0$---this follows from:
$$
\frac{k+1}{\root4\of{1+\sqrt2} (k+\tfrac12)}
< \frac{k+1}{k+\tfrac12} \le 2 < e^{3\pi/8} \le e^{\tfrac\pi2 x(k+\tfrac34)}
= e^{\tfrac\pi2 x[(k+1)^2-(k+\tfrac12)^2]}.
$$
This proves (6.6).
Now as $\vartheta_3(0,2ix)+(\sqrt2-1)\vartheta_2(0,2ix)$ is greater than 1 (since
$\vartheta_3(0,2ix)>1$), we can multiply it for free on the right side of (6.6) so that 
this product becomes 
$$
\multline
\left( \vartheta_3(0,2ix)-(1+\sqrt2)\vartheta_2(0,2ix) \right)
\left( \vartheta_3(0,2ix)+(\sqrt2-1)\vartheta_2(0,2ix) \right)
\\
\ = \ \left( \vartheta_3(0,2ix)-\vartheta_2(0,2ix) \right)^2 - 2\vartheta_2(0,2ix)^2
\ = \ \vartheta_3(\tfrac\pi2,\tfrac i2 x)^2 - 2\vartheta_2(0,2ix)^2
\endmultline
$$
using the identity $\vartheta_3(w,u) = \vartheta_3(2w,4u) + \vartheta_2(2w,4u)$ and the
fact that $\vartheta_2(\pi,v)=-\vartheta_2(0,v)$. Therefore, (6.6) yields
$$
h(x) < \vartheta_3(\tfrac\pi2,\tfrac i2 x)^2 - 2\vartheta_2(0,2ix)^2
$$
which is exactly (6.5).  \qed

\enddemo

\bigpagebreak

\proclaim{Lemma 6.3} \ $\vartheta_3(0,2i) = (1+\sqrt2) \vartheta_2(0,2i)$.
\endproclaim
\demo{Proof}
Let $x = \vartheta_3(0,2i)/ \vartheta_2(0,2i)$.  Since $x>0$, it will suffice to show
that $x$ satisfies the equation $(x-1)^2=2$, or $x(x-1)=x+1$.  This becomes
$$
\vartheta_3(0,2i)\bigl(\vartheta_3(0,2i)-\vartheta_2(0,2i)\bigr) \ = \ 
\vartheta_2(0,2i)\bigl(\vartheta_3(0,2i)+\vartheta_2(0,2i)\bigr).
\tag6.7
$$
Using the Theta function identities
$$
\vartheta_3(w,u) = \vartheta_3(2w,4u) + \vartheta_2(2w,4u), \qquad
\vartheta_4(v,u) =  \vartheta_3(2v,4u) - \vartheta_2(2v,4u)
$$
(6.7) becomes
$$
\vartheta_3(0,2i)\,\vartheta_4(0,\tfrac i2) \ = \ 
\vartheta_2(0,2i)\,\vartheta_3(0,\tfrac i2).
\tag6.8
$$
Applying the following inversion formulas
$$
\vartheta_4(z,t) = (-it)^{-1/2}\, e^{z^2/(\pi it)}\,
\vartheta_2\left(\tfrac z t, - \tfrac 1 t\right), \qquad
\vartheta_3(z,t) = (-it)^{-1/2}\, e^{z^2/(\pi it)}\,
\vartheta_3\left(\tfrac z t, - \tfrac 1 t\right),
$$
one obtains $\vartheta_4(0,\tfrac i2)=\sqrt2 \vartheta_2(0,2i)$ and 
$\vartheta_3(0,\tfrac i2)=\sqrt2 \vartheta_3(0,2i)$.  Now it is clear that (6.8) holds.\qed
\enddemo

\remark{Remark} From the proof, and using $\vartheta_3(z+\frac\pi2,u) = \vartheta_4(z,u)$,
one also has
$$
\frac{\vartheta_3(0,\tfrac i2)}{\vartheta_3(\tfrac\pi2,\tfrac i2)}
= \frac{\sqrt2 \vartheta_3(0,2i)}{\vartheta_4(0,\tfrac i2)}
= \frac{\sqrt2 \vartheta_3(0,2i)}{\sqrt2 \vartheta_2(0,2i)}
= 1+\sqrt2.
$$
\endremark

We have thus proved the following result.

\proclaim{Theorem 6.4} Let $\rho>1$ and $U,V$ be canonical unitaries generating the 
rotation algebra $A_\rho$ with $VU=e(\rho)UV$.  Then the element
$$
\sum_{m,n}  e^{-\tfrac\pi2\rho(m^2+n^2)} e(\tfrac12\rho mn)  U^n V^m
$$
is positive and invertible in $A_\rho$ (and invariant under the Fourier transform).
\endproclaim

We remark that for $0<\rho<1$ the element in Theorem 6.4 cannot be invertible, for then
it would mean that by using an appropriate Rieffel framework one could construct a 
projection whose (normalized) trace is greater than 1.  For $\rho=1$ the element is 
a function in $C(\Bbb T^2)$ and one can check that it is zero when $U$ and $V$ 
are replaced by $-1$.

\newpage

\subhead \S7. APPROXIMATE CENTRALITY \endsubhead

We will now prove that the Fourier invariant projection $e=\Dinner \xi\xi$, where
$\xi = fb$, whose existence is now assured by Section 6, is approximately central.

\proclaim{Proposition 7.1} Let $\theta$ be any irrational number and let $\{p/q\}$
be a sequence of positive rationals with $q\to\infty$, such that 
$q|q\theta-p|< \kappa < 1$ for some $\kappa<1$.  (For example, this is satisfied with 
$\kappa=1/2$ by at least one of any two consecutive convergents of $\theta$.)
Then the associated Fourier invariant projection $e$ (of trace $q|q\theta-p|$) is 
approximately central.
\endproclaim

Before proving this let us first prove that the primitive form $X:=\Dinner ff$ of the 
projection is itself approximately central. From the end of Section 4 we have
$$
X = \frac1{\beta^2} \sum_{m,n} e(\tfrac12\alpha^2q^2mn) e^{-\tfrac\pi2\alpha^2q^2(m^2+n^2)} 
\mu_{qm,qn} U_2^{qn} U_1^{qm}.
$$
Since $X$ is Fourier invariant, it is enough to show that $\|U_1XU_1^*-X\| \to 0$ 
(inserting $\sigma$ inside the norm gives $\|U_2XU_2^*-X\| \to 0$).  Now
$$
U_1XU_1^* - X = \frac1{\beta^2} 
\sum_{m,n} e(\tfrac12\alpha^2q^2mn) e^{-\tfrac\pi2 \alpha^2q^2(m^2+n^2)} 
\mu_{qm,qn} \, (\lambda^{qn}-1) \, U_2^{qn} U_1^{qm}
$$
hence
$$
\|U_1XU_1^* - X\| \le \frac1{\beta^2}
\sum_{m,n} e^{-\tfrac\pi2 \alpha^2q^2(m^2+n^2)} \, |\lambda^{qn}-1| = 
\frac2{\beta^2} \vartheta_3(0,\tfrac i2\alpha^2q^2) 
\sum_{n=1}^\infty e^{-\tfrac\pi2 \alpha^2q^2n^2} \, |\lambda^{qn}-1|.
$$
We have $\lambda^{qn}=e(\alpha^2qn)$, and as $|e(t)-1|\le 2\pi t$ for $t\ge0$, one gets
(since $\beta=\tfrac1{q\alpha}$)
$$
\sum_{n=1}^\infty e^{-\tfrac\pi2 \alpha^2q^2n^2} \, |\lambda^{qn}-1| 
\le
2\pi \alpha^2q \sum_{n=1}^\infty n e^{-\tfrac\pi2 \alpha^2q^2n^2}
=
\frac{2\pi}{q\beta^2} \sum_{n=1}^\infty n e^{-\tfrac\pi{2\beta^2}n^2}.
$$
It is not hard to check that one has 
$$
\sum_{n=1}^\infty n e^{-\tfrac\pi{2\beta^2} n^2} \ \le \ 
\frac1\pi \beta^2 + \frac2{\sqrt{\pi e}} \beta
$$
(e.g., see Lemma A.1 of [\SWd]). 
Using the inequality $\vartheta_3(0,ix) \le 1 + \tfrac1{\sqrt x}$, one gets (as $\beta>1$)
$$
\|U_1XU_1^* - X\| \ \le \
\frac{4\pi}{q\beta^2} (1+\beta\sqrt2) \left(\frac1\pi + \frac2{\beta\sqrt{\pi e}} \right)
< \frac{12\pi}q.
$$
Therefore, $X$ is approximately central.  Next we need the following.

\proclaim{Lemma 7.2} Let $\theta$ be any irrational number.  For any sequence of 
positive rationals $\{p/q\}$ satisfying $q|q\theta-p|< \kappa < 1$, the norms of 
$\Dinner ff$ and $\Dinner ff^{-1}$ are bounded above by a constant that is independent
of $p/q$ (but which depends on $\kappa$).
\endproclaim
\demo{Proof} 
Since $\mu_\xi: eA_\theta e \to C^*(D^\perp)$ ($\mu_\xi(x)=\Dpinner \xi{x\xi}$) is an 
isomorphism of C*-algebras, and $\mu_\xi(\Dinner ff) = \Dpinner ff$, it is enough 
to show that the norms of $\Dpinner ff$ and $\Dpinner ff^{-1}$ are bounded above.
From (6.4) we obtained
$$
\| R^{-1} \rho_0^{-1} \Dpinner ff  - 1 \| \le \Cal E(\beta^2) 
$$
where
$$ 
\Cal E(\beta^2) = \frac{8\pi (\beta^2-1) \Psi(\tfrac12\beta^2) \Psi(2\beta^2) +
2\vartheta_2(0,2i\beta^2)^2}
{\vartheta_3(\tfrac\pi2,\tfrac i2\beta^2)^2}\ < \ 1
$$
holds for any $\beta^2 = \tfrac1{q|q\theta-p|}>1$. Since $\Cal E(\beta^2)$ is a continuous 
function, and since it goes to zero as $\beta\to\infty$, its maximum over the 
$\beta^2$-interval $[\tfrac1\kappa,\infty)$ is some constant $C<1$ independent of the
sequence $p/q$ considered. Hence for all $p/q$ such that $q|q\theta-p|< \kappa$ one has
$$
\| R^{-1} \rho_0^{-1} \Dpinner ff  - 1 \| \le C
\tag7.1
$$
from which one has 
$$
\Dpinner ff^{-1}\rho_0R = \sum_{n=0}^\infty (1-R^{-1} \rho_0^{-1} \Dpinner ff)^n
$$
so that
$$
\|\Dpinner ff^{-1}\| \ 
\le \ \|\rho_0^{-1}\|^2 \sum_{n=0}^\infty C^n \ = \ \frac{\|\rho_0^{-1}\|^2}{1-C}
$$
In addition, we have $\|R^{-1}\| = \|\rho_0^{-1}\| = 
\vartheta_3(\tfrac\pi2,\tfrac i2\beta^2)^{-1} > 0$, by (6.1), is also bounded above 
by a constant independent of $p/q$ (in view of our hypothesis on the rationals $p/q$).
(Note that $x\mapsto \vartheta_3(\tfrac\pi2,\tfrac i2x)$ is an increasing function over
$x\ge1$ and bounded above by 1.)
Using (7.1) one easily obtains an upper bound for $\|\Dpinner ff\|$. \qed
\enddemo

\bigpagebreak

\demo{Proof of Proposition 7.1} Fix $\epsilon>0$.
With $X=\Dinner ff$, one checks by induction that
$$
X^k = \Dinner{\xi c^k}\xi
$$
for each integer $k\ge1$, where $c:=b^{-2}:=\Dpinner ff$.  Therefore for any polynomial
$P(x)$ such that $P(0)=0$ one has
$$
P(X) = \Dinner{\xi P(c)}\xi.
$$
By Lemma 7.2 there is a closed interval $[r,s]$, where $r>0$, that
contains the spectrum of $X$ (and hence also that of $U_1XU_1^*$) for all $p/q$ such that
$q|q\theta-p|< \kappa$.  This allows one to obtain a sequence of polynomials $\{P_n(x)\}$
such that $P_n(0)=0$ for each $n$ and $P_n(x) \to 1$ uniformly on $[r,s]$.  Since $X$ is
supported on the projection $e$, and invertible in $eA_\theta e$, one gets 
$\|P_n(X)-e\| = \|P_n(X)-1(X)\| \le \|P_n-1\|_{[r,s]} < \epsilon$ for large enough $n$.
Since we showed that $\|U_1XU_1^*-X\| \to 0$ as $q\to\infty$, we can find $q$ large 
enough (and depending on $n$) so that
$$
\|P_n(U_1XU_1^*) - P_n(X)\| < \epsilon.
$$
Therefore, for large $q$ one has
$$
\|U_1eU_1^* - e \|  \le  \|U_1eU_1^* - U_1P_n(X)U_1^*\| + \|P_n(U_1XU_1^*) - P_n(X)\|
+ \|P_n(X) - e\|  < 3\epsilon.
$$
Since $e$ is Fourier invariant this implies that $\|U_2eU_2^* - e \|< 3\epsilon$.  \qed
\enddemo

\newpage

\subhead \S8. THE CUT DOWN APPROXIMATION \endsubhead

We need only approximate the cut down ``$eU_1e$'' $:= \Dpinner \xi{U_1\xi}$ by elements
of the $q\times q$ matrix algebra generated by $V_3,V_4$, since by taking the Fourier 
transform of the result one gets the approximation for ``$eU_2e$'' in the same matrix 
algebra (as $e$ is Fourier invariant).  From (5.8) we obtained
$$
\Dpinner f{U_1f} = \mu_0^{1/2} e(\tfrac K{2q})  V_4^{c_4}  V_3^{c_3}
\sum_{m,n}  e^{-\tfrac\pi2[(\alpha+\beta m)^2+\beta^2n^2]} 
e(\tfrac12\beta^2 mn) e(-\tfrac n{2q}) W_2^n W_1^m
$$
where, 
$$
c_3 = \Delta'(-r_1s_3+r_2s_4), \qquad c_4 = \Delta'(-r_4s_4+r_3s_3)
$$
and $r_j,s_j$ are given by (3.3).  Write
$$
\Dpinner f{U_1f} = \mu_0^{1/2} e(\tfrac K{2q})  V_4^{c_4}  V_3^{c_3} B
$$
where
$$
B := \sum_{m,n}  e^{-\tfrac\pi2[(\alpha+\beta m)^2+\beta^2n^2]} 
e(\tfrac12\beta^2 mn) e(-\tfrac n{2q}) W_2^n W_1^m
$$
Since $\xi := fb$, we have
$$
\Dpinner \xi{U_1\xi} = \mu_0^{1/2} e(\tfrac K{2q})  bV_4^{c_4}  V_3^{c_3} \,Bb.
$$
To show that $\Dpinner \xi{U_1\xi}$ is approximately equal to 
$\mu_0^{1/2} e(\tfrac K{2q}) V_4^{c_4}  V_3^{c_3}$, which belongs to the Fourier invariant 
matrix algebra $C^*(V_3,V_4)$, it will suffice to show the following:

\medpagebreak

\itemitem{1)} $\|B-b^{-2}\| \to 0$,
\itemitem{2)} $\|V_4^{c_4}V_3^{c_3}b V_3^{-c_3}V_4^{-c_4} - b \| \to 0$,

\medpagebreak

\noindent 
as $q\to\infty$ through a sequence of rationals $p/q$ satisfying the hypothesis of 
Proposition 7.1-----and for which the norms of $b$ and $b^{-1}$ are uniformly bounded, 
in view of Lemma 7.2.

\demo{Proof of 1)} Fix $\epsilon>0$ and let $N$ be any positive integer. We have
$$
\align
B - b^{-2} &= 
\sum_{m,n}  e^{-\tfrac\pi2\beta^2(m^2+n^2)} 
e(\tfrac12\beta^2 mn) \left( e^{-\tfrac\pi2(\alpha^2+\tfrac2qm)} e(-\tfrac n{2q} ) - 1 
\right) W_2^n W_1^m
\\
&= \sum_{|m|\le N, |n|\le N} + \sum_{|m|>N, |n|\le N} + \sum_{|n|>N, |m|\le N} + \sum_{|m|>N, |n|> N}.
\endalign
$$
For the second sum we have (as $\beta>1$)
$$
\align
\left\|  \sum_{|m|>N, |n|\le N} \right\|
&\le 
\sum_{|m|>N, |n|\le N} e^{-\tfrac\pi2\beta^2(m^2+n^2)} 
\left( e^{-\tfrac\pi2(\alpha^2+\tfrac2qm)}  + 1 \right)
\\
&< \sum_{|m|>N} e^{-\tfrac\pi2 m^2} 
\left( e^{-\tfrac\pi2(\alpha^2+\tfrac2qm)}  + 1 \right)
 \sum_{|n|\le N} e^{-\tfrac\pi2 n^2} 
\\
&< \vartheta_3(0,\tfrac i2)  \sum_{|m|>N} e^{-\tfrac\pi2 m^2} 
\left( e^{-\tfrac\pi2(\alpha^2+\tfrac2qm)}  + 1 \right)
\endalign
$$
which is clearly less than $\epsilon$ for $N>N_1$ for large enough $N_1$ independent of
$q$ (and $\alpha \to0$ as $q\to\infty$).  This follows from the following inequality
$$
\sum_{|m|>N} e^{-\tfrac\pi2 m^2} e^{-\tfrac\pi qm} \ < \ 2 
\sum_{k=N}^\infty e^{-\tfrac\pi2 k^2}
$$
which holds for all $q\ge1$. Similarly, for the third sum
$$
\align
\left\|  \sum_{|n|>N, |m|\le N} \right\|
&\le 
\sum_{|n|>N, |m|\le N} e^{-\tfrac\pi2\beta^2(m^2+n^2)} 
\left( e^{-\tfrac\pi2(\alpha^2+\tfrac2qm)}  + 1 \right)
\\
&< \sum_{|m|\le N} e^{-\tfrac\pi2 m^2} 
\left( e^{-\tfrac\pi2(\alpha^2+\tfrac2qm)}  + 1 \right)
 \sum_{|n|>N} e^{-\tfrac\pi2 n^2} 
\endalign
$$
the sum over $m$ here is bounded by a positive number independent of $q$ and the sum
over $n$ goes to zero, so that the preceeding norm is less than  $\epsilon$ for $N>N_2$
for some $N_2$ independent of $q$.  The fourth sum similarly is less than $\epsilon$ 
for $N>N_3$ for some $N_3$ independent of $q$.  Letting $N = \max(N_1,N_2,N_3)$ one has
$$
\left\|  \sum_{|m|,|n|\le N} \right\| <
\sum_{|m|,|n|\le N} e^{-\tfrac\pi2(m^2+n^2)} 
\left| e^{-\tfrac\pi2(\alpha^2+\tfrac2qm)} e(-\tfrac n{2q} ) - 1 \right|
$$
and this sum (now $N$ being fixed and independent of $q$) can be made less than 
$\epsilon$ by taking $q$ large enough.  This proves 1).  \qed
\enddemo

\medpagebreak

\demo{Proof of 2)} First, one easily checks that
$$
V_4^{c_4}  V_3^{c_3} W_j V_3^{-c_3}  V_4^{-c_4} = \nu_j W_j
$$
where $\nu_j := e(\tfrac pq(c_3b_j-c_4a_j))$.  Since from (3.5) (where $u_3=s_3,\ u_4=s_4$)
we have $a_1=-cc_3$ and $b_1=-cc_4$, one has $c_3b_1-c_4a_1=0$ so that $\nu_1=1$.  Also, as in 
the calculation at the end of the proof of Lemma 5.1
$$
\align
c_4a_2-c_3b_2 &=
c{\Delta'}^2(-r_4s_4+r_3s_3)(r_1s_1-r_2s_2) - c{\Delta'}^2(-r_1s_3+r_2s_4)(r_4s_2-r_3s_1)
\\
&= c{\Delta'}^2 
\Bigl( (r_1s_3-r_2s_4)(r_4s_2-r_3s_1) - (r_1s_1-r_2s_2)(r_4s_4-r_3s_3) \Bigr)
\\
&= -c{\Delta'}^2 (r_1r_4-r_2r_3)(s_1s_4-s_2s_3)
\\
&= -c{\Delta'}^2 \Delta^2
\\
&\equiv -c \mod q
\endalign
$$
Hence $\nu_2 = e(\tfrac pq(c_3b_2-c_4a_2)) = e(\tfrac{pc}q) = e(\tfrac1q)$, since 
$pc\equiv1$.  Therefore, we have
$$
V_4^{c_4}  V_3^{c_3} b^{-2} V_3^{-c_3}  V_4^{-c_4} - b^{-2} =
\sum_{m,n}  e^{-\tfrac\pi2\beta^2(m^2+n^2)} e(\tfrac12\beta^2 mn) (e(\tfrac nq)-1) 
W_2^n W_1^m
$$
hence
$$
\|V_4^{c_4}  V_3^{c_3} b^{-2} V_3^{-c_3}  V_4^{-c_4} - b^{-2}\|
\le 
\sum_{m,n}  e^{-\tfrac\pi2(m^2+n^2)} |e(\tfrac nq)-1|
\le
\frac{2\pi}q \sum_{m,n} |n|e^{-\tfrac\pi2(m^2+n^2)}
$$
which goes to 0 as $q\to\infty$.  Since the norms of $b^2$ and $b^{-2}$ are
uniformly bounded by Lemma 7.2 (for rationals $p/q$ satisfying the hypothesis therein),
one has thereby proved 2).  \qed
\enddemo
We thus conclude that
$$
\Dpinner \xi{U_1\xi} \approx \mu_0^{1/2} e(\tfrac K{2q})V_4^{c_4}V_3^{c_3}.
$$
for large enough $q$.
Taking the Fourier transform ($\sigma'$) of this approximation and using $\ft\xi = \xi W_0$
from Section 4, one obtains
$$
W_0^*  \Dpinner {\xi}{U_2\xi} W_0 =  \Dpinner {\xi W_0}{U_2\xi W_0} 
= \Dpinner {\ft\xi}{U_2\ft\xi} = \sigma'(\Dpinner {\xi}{U_1\xi})
\approx \mu_0^{1/2} e(\tfrac K{2q}) V_3^{-c_4} V_4^{c_3}
$$
and thus
$$
\text{``}eU_2e\text{''} := \Dpinner {\xi}{U_2\xi} 
\approx \mu_0^{1/2} e(\tfrac K{2q}) W_0 V_3^{-c_4} V_4^{c_3} W_0^* 
$$
and the latter element is in the matrix algebra $C^*(V_3,V_4) \cong M_q(\Bbb C)$ 
since it contains $W_0$.  This completes the proof of (4) of Theorem 1.1.

\newpage

\subhead \S9. THE AF STRUCTURE OF $A_\theta \rtimes \Bbb Z_4$ \endsubhead

In this section we prove that for a dense $G_\delta$ set $\Cal G$ of irrationals 
$\theta$ (to be described shortly) the C*-algebra $A_\theta \rtimes\Bbb Z_4$ is 
an AF-algebra.  To do this we will first need to show that it is tracially AF and 
satisfies UCT.

\proclaim{Theorem 9.1}
For $\theta$ in the dense $G_\delta$ set $\Cal G$, the C*-algebra 
$A_\theta \rtimes\Bbb Z_4$ is tracially AF.  In fact, for each $\epsilon>0$ there is
a Fourier invariant projection $e$ in $A_\theta$ such that $\vartau(1-e)<\epsilon$,
$e$ is approximately central, and $eU_1e, eU_2e$ are approximately contained in a Fourier
invariant $q\times q$ matrix algebra that has $e$ as its unit. 
\endproclaim

Using H.~Lin's Theorems 3.4 and 3.6 of [\HLb], Theorem 9.1 yields 

\proclaim{Theorem 9.2} 
For $\theta \in \Cal G$, the C*-algebra 
$A_\theta \rtimes\Bbb Z_4$ has real rank zero, stable rank one, is quasidiagonal, 
and its $K_0$ group is weakly unperforated.
\endproclaim

Combining this together with Theorem 9.4 below, which shows that $A_\theta \rtimes\Bbb Z_4$
satisfies UCT, we can use a recent result of H.~Lin's-----namely that every unital 
separable simple nuclear tracially AF algebra that satisfies UCT is isomorphic to an 
AH-algebra with slow dimension growth (see Remark 4.6 of [\HLa])-----to deduce that
$A_\theta \rtimes\Bbb Z_4$ is an AH-algebra with slow dimension growth (for $\theta \in
\Cal G$).
This puts the algebra $A_\theta \rtimes\Bbb Z_4$ into the classification class of Elliott 
and Gong (now that it has real rank zero and stable rank one).  
In earlier work [\SWb] we showed that its $K_1$ is zero and $K_0 = \Bbb Z^9$ 
for $\theta$ in another dense $G_\delta$. Since the crossed product is simple, unital, stably
finite, has weakly unperforated $K_0$ (for a dense $G_\delta$), and has a unique trace 
state $\vartau$, by a theorem of Blackadar and Handelman (see [\BB], 6.9.2) for two 
positive 
classes $x,y$ in $K_0$ one has $x<y$ iff $\vartau(x)<\vartau(y)$.  This shows that the 
Riesz interpolation property is automatically satisfied, and since $K_0$ is torsion 
free and weakly unperforated it must be unperforated, therefore by a theorem of 
Effros-Handelman-Shen ([\BB], 7.4.1) the Elliott invariant of 
$A_\theta \rtimes\Bbb Z_4$ (now just consisting of the scaled ordered group structure 
of $K_0$) is a dimension group, and so is equal to that of an AF-algebra. 
Therefore, thanks to the Elliott-Gong Classification Theorem [\EG] one obtains:

\proclaim{Theorem 9.3} 
For a dense $G_\delta$ set of irrationals $\theta$, the C*-algebra 
$A_\theta \rtimes\Bbb Z_4$ is an AF-algebra.
\endproclaim

\medpagebreak

It only remains to prove Theorem 9.1 and show that $A_\theta \rtimes\Bbb Z_4$ satisfies UCT.
First, however, let us settle the UCT property.

\proclaim{Theorem 9.4}
For each $\theta$ in $(0,1)$ the C*-algebra $A_\theta \rtimes\Bbb Z_4$ satisfies UCT.
\endproclaim

\demo{Proof}
Let $G$ denote the semidirect product of the discrete Heisenberg group, generated by
elements $u,v$ such that $vu=cuv$, where $c$ is central, by the canonical order four 
automorphism (the ``Fourier transform'') $u\mapsto v \mapsto u^{-1}$.  
Since $G$ is an amenable group, we can use Tu's Theorem [\TU]-----which says that the 
group C*-algebra $C^*(G)$ of an amenable groupoid belongs to the so-called bootstrap
category ([\BB], 22.3.4)-----to deduce that $C^*(G)$ satisfies UCT.
The algebra $C^*(G)$ can be viewed as a continuous field of C*-algebras over its 
primitive spectrum ([\JD], 10.5.1), which is homeomorphic to the unit circle, and 
which we view as the closed interval $[0,1]$ with its endpoints identified.
The evaluation map on $C^*(G)$ at $\theta\in(0,1)$ gives a natural surjection onto 
the fiber $A_\theta \rtimes\Bbb Z_4$,
whose kernel is the ideal $J_\theta$ of all continuous sections vanishing at $\theta$.
Note that if $\theta$ is rational then $J_\theta$ satisfies UCT since in this case 
the crossed product $A_\theta \rtimes\Bbb Z_4$ is of type I.
By the ``two out of three'' theorem, it is enough to show $J_\theta$ satisfies UCT. 
Consider the following two steps.

\itemitem{(1)} Consider the ideal $J_I$ of all sections in $C^*(G)$ that vanish on an
open subinterval $I$ of $(0,1)$ with rational endpoints.  This ideal is a direct
summand of the ideal of all sections that vanish at its endpoints (which
satisfies UCT), hence $J_I$ satisfies UCT (see [\BB], end of 23.10).  
If one takes a continuous section that  
vanishes at two rational points $r < s$, say, then it can be restricted to
$[r,s]$ making it zero outside this interval, and it can also be restriced
to $[0,r] \cup [s,1]$ (with $0$ and $1$ identified) and making it zero on $[r,s]$.
These are continuous sections by the local property for continuous
fields ([\JD], 10.1.2 (iv)), and the original section is their (direct) sum.

\itemitem{(2)} 
If a continuous section $\zeta$ is zero at $\theta$, then it is close in norm to a
section that vanishes in some small enough open interval (with rational
endpoints) that contains $\theta$.  This can be done because the continuous
sections of any continuous field of C*-algebras is closed under
multiplication by scalar continuous functions ([\JD], 10.1.9).  So if we take
a scalar function $g$ that is $1$ off some small open interval containing $\theta$,
but which is zero on a smaller subinterval containing $\theta$, then $g\zeta$ is close
to $\zeta$ in norm and $g\zeta$ is in $J_I$ for suitable $I$.

\noindent By (1) and (2), it follows that given any finite number of continuous
sections in $J_\theta$, then each is close in norm to a section that belongs to
some $J_I$ for some small enough open interval (with rational endpoints) containing 
$\theta$.  Since $J_I$ satisfies UCT, and is a C*-subalgebra of $J_\theta$, one applies 
Dadarlat's recent theorem [\MD] (Theorem 1.1) to deduce that $J_\theta$ satisfies UCT.
(In his paper, Dadarlat proves that if a nuclear separable C*-algebra can be 
approximated by C*-subalgebras satisfying UCT, then it also satisfies UCT.)  \qed
\enddemo

\remark{Remark}
The author is thankful to Chris Phillips for suggesting the idea of using Tu's Theorem 
and viewing $C^*(G)$ as a continuous field of C*-algebras.
\endremark

\remark{Remark}
The same proof shows that the crossed products $A_\theta \rtimes_\alpha \Bbb Z_6$ and 
$A_\theta \rtimes_{\alpha^2}\Bbb Z_3$, where $\alpha$ is the canonical order 6 automorphism 
given by $\alpha(U)=V,\ \alpha(V)=U^{-1}V$, satisfy the Universal Coefficient Theorem
for any $\theta$, by taking $G$ to be the appropriate semi-direct product of the discrete 
Heisenberg group by the canonical action of $\Bbb Z_6$, respectively $\Bbb Z_3$-----which 
is amenable.
\endremark

\bigpagebreak

It now remains to prove Theorem 9.1.  We specify the dense $G_\delta$ set $\Cal G$
of irrationals as follows.  Pick any two sequences of positive
integers $\{M_k\}, \{N_k\}$ that go to infinity.  Look at the set $\Cal F$ of all 
irrationals $\theta$ whose continued fraction expansion $[a_0,a_1,a_2,\dots]$ is 
such that for each $k$ the triple $N_k,1,M_k$ appears infinitely often in the 
sequence $[a_0,a_1,a_2,\dots]$.  (Thus, for each $k$ there
are infinitely many $n$ with $a_n=N_k, a_{n+1}=1, a_{n+2}=M_k$.)  We define the class
$\Cal G$ to be the union of the classes $\Cal F$ over all pairs of sequences 
$\{M_k\}, \{N_k\}$.

\smallpagebreak

For fixed sequences $\{M_k\}, \{N_k\}$, that the aforementioned set $\Cal F$ is 
a dense $G_\delta$ can be seen as follows.  First, for any three positive integers 
$a,b,c$ the set $S$ of irrationals $\theta$ such that the triple $a,b,c$ appears 
consecutively in its continued fraction expansion $[a_0,a_1,a_2,\dots]$ is a dense 
$G_\delta$. Indeed, given such $\theta$, let us say $a_n=a,a_{n+1}=b,a_{n+2}=c$, 
it is clear that there is $\epsilon>0$ such that each 
irrational $\theta'$ in $(\theta-\epsilon,\theta+\epsilon)$ has the same continued fraction
coefficients as $\theta$ up to $a_{n+2}$, and so in particular $a,b,c$ occurs in the
continued fraction expansion of $\theta'$.  Hence the set $S$ is open and dense in the 
irrationals, so it is a dense $G_\delta$.  Since the countable intersection of dense 
$G_\delta$'s is again a dense $G_\delta$, the set of irrationals such that the triple 
$a,b,c$ appears infinitely often in its continued fraction is still a dense $G_\delta$. 
This shows that the class $\Cal G$ is a dense $G_\delta$.

\demo{Proof of Theorem 9.1}  Fix $\theta\in \Cal G$.
From the theory of continued fractions it is well-known that each convergent 
$p_n/q_n$ of $\theta$ satisfies
$$
q_n\theta-p_n = \frac{(-1)^n}{\xi_{n+1}q_n + q_{n-1}} 
$$
where $\xi_{n+1} = a_{n+1} + \tfrac1{a_{n+2}+\dots} = [a_{n+1},a_{n+2},\dots]$.
Write the latter as $\xi_{n+1} = a_{n+1} + \tfrac1{a_{n+2}+\eta}$ where $0<\eta<1$.
For $n$ such that $a_n=N_k,\, a_{n+1}=1,\, a_{n+2}=M_k$, we have 
$$
\beta^2 := \frac1{q_n|q_n\theta-p_n|} = \xi_{n+1} + \frac{q_{n-1}}{q_n} > \xi_{n+1} =
a_{n+1} + \frac1{a_{n+2}+\eta} = 1 + \frac1{M_k+\eta} >  1 + \frac1{M_k+1}.
$$
On the other hand, since 
$$
\xi_{n+1} = 1 + \frac1{M_k+\eta} < 1 + \frac1{M_k}, \qquad 
\frac{q_{n-1}}{q_n} = \frac{q_{n-1}}{a_nq_{n-1}+q_{n-2}}
= \frac1{a_n+\tfrac{q_{n-2}}{q_{n-1}}} < \frac1{N_k}
$$
we also have
$$
\beta^2 = \xi_{n+1} + \frac{q_{n-1}}{q_n} < 1 + \frac1{M_k} + \frac1{N_k}.
$$
Therefore, for each $k$ there are infinitely many convergents $p/q$ ($=p_n/q_n$) such that
$$
1 + \frac1{M_k+1} \ < \ \beta^2 = \frac1{q|q\theta-p|} \ <\  1 + \frac1{M_k} + \frac1{N_k}.
\tag*
$$
The important thing we get here is that $\beta$ can be made as close to 1 as we want
(by choosing $k$ large enough) but at the same time we can keep it away from 1 for 
infinitely many $q$. Fix $\epsilon>0$, and fix an integer $k$ such that
$$
\vartau(1-e) = 1 - q|q\theta-p| \ < \ 1 - \frac1{1 + \frac1{M_k} + \frac1{N_k}} \ 
< \ \epsilon
$$
which now holds for infinitely many convergents $p/q$ independent of $k$, by (*).  
The hypothesis of Theorem 7.1 is satisfied by (*) so that the Fourier invariant 
projection $e$ is approximately central.  In Section 8 it  was already shown that $e$ is a 
point projection, in that the corner algebra $eA_\theta e$ can be approximated by a 
Fourier invariant $q\times q$ matrix subalgebra (whose unit is also $e$).  This 
proves Theorem 9.1. \qed
\enddemo

\remark{Remark}
To prove Theorem 9.1 for all irrational $\theta$, one would need at least two point 
projections the sum of whose traces is close to 1.  This seems quite conceivable to do
(but this paper is long enough!).  Using the ideas in this paper, if one can prove that
$A_\theta \rtimes \Bbb Z_4$ is tracially AF, $K_1=0$, and $K_0=\Bbb Z^9$ for all irrational
$\theta$, then it follows that $A_\theta \rtimes \Bbb Z_4$ is an AF-algebra for all
irrational $\theta$.  This we leave for future work.
\endremark


				\enddocument